\documentclass[letter,11pt]{article}

\usepackage[margin=1in]{geometry}
\usepackage[useregional]{datetime2}
\usepackage[english]{babel}
\usepackage[utf8]{inputenc}
\usepackage[T1]{fontenc}
\usepackage{amsfonts}
\usepackage{amsmath}
\usepackage{amssymb}
\usepackage{amsthm}
\usepackage{mathrsfs}
\usepackage{dsfont}
\usepackage{cite}
\usepackage{amsxtra,amscd}
\usepackage{hyperref}
\usepackage{eurosym}
\usepackage{multicol}
\usepackage{xcolor}
\usepackage{ulem}
\usepackage{bbm}
\usepackage{csquotes}
\usepackage{extarrows}

\usepackage{xcolor, graphicx,changepage,xcolor, float, array, amssymb, amsmath, mathrsfs, amsfonts, hyperref, url, amsthm, tikz-cd, fancyhdr, enumerate, bbm}
\makeatletter
\newcommand*{\rom}[1]{\expandafter\@slowromancap\romannumeral #1@}
\makeatother

\makeatletter
\renewcommand*\env@matrix[1][*\c@MaxMatrixCols c]{%
	\hskip -\arraycolsep
	\let\@ifnextchar\new@ifnextchar
	\array{#1}}
\makeatother

\newcommand{\C}{\mathbb{C}}
\newcommand{\N}{\mathbb{N}}

\newcommand{\Z}{\mathbb{Z}}

\newcommand{\uproman}[1]{\uppercase\expandafter{\romannumeral#1}}

\theoremstyle{definition}
\newtheorem{theorem}{Theorem}[section]
\newtheorem{lemma}[theorem]{Lemma}
\newtheorem{definition}[theorem]{Definition}

\newtheorem{corollary}[theorem]{Corollary}
\newtheorem{proposition}[theorem]{Proposition}
\newtheorem{example}[theorem]{Example}

\theoremstyle{remark}
\newtheorem*{remark}{Remark}

\usepackage{setspace}
\usepackage{ulem} 

\title{On the Jacobs-de Leeuw-Glicksberg decomposition for representations of semigroups}
\author{Micky Barthmann\footnote{Chemnitz University of Technology, Faculty of Mathematics,  Reichenhainer Straße 41, 09126 Chemnitz, Germany. Email: micky.barthmann@mathematik.tu-chemnitz.de},\ \ Sohail Farhangi\footnote{University of Adam Mickiewicz, Department of Mathematics and Informatics, ulica Wieniawskiego 1, 61-712 Poznań, Poland. Email: sohail.farhangi@gmail.com}\ \footnote{Beijing Institute of Mathematical Sciences and Applications, Beijing, P.R.C. 101408.},\ \ and
Yulia Kuznetsova\footnote{Universit\'e Marie et Louis Pasteur, CNRS, LmB (UMR 6623), F-25000 Besançon, France}}
\date{\today}

\begin{document}

\maketitle
\begin{abstract}
   We develop a systematic study of the Jacobs--de Leeuw--Glicksberg decomposition for JdLG-admissible semigroup representations on Banach spaces, with emphasis on its structural, ergodic, and combinatorial aspects.
    For a JdLG-admissible representation $\pi$ of a semigroup $S$ on a Banach space $E$, we show that the reversible part is weakly equivalent to a unitary representation on a Hilbert space that decomposes as a direct sum of finite-dimensional representations.
    We also characterize the almost weakly stable part in terms of the unique invariant mean on the space of weakly almost periodic functions. 
    When $S$ is a bi-amenable measured semigroup, we obtain further characterizations of the almost weakly stable part using invariant means and averages along F\o lner sequences. 
    In addition, we describe, in terms of ultrafilters, the unique projection onto the reversible part whose kernel is the almost weakly stable part, and derive several combinatorial consequences.
    A range of examples are given to demonstrate the sharpness of our results.
    Many of our theorems extend known features of the compact--weak mixing decomposition for unitary operators; although some auxiliary results were previously known or part of the folklore, their extension to this general setting, particularly for nonseparable Banach spaces, requires new ideas.

\end{abstract}
\section{Introduction}
The compact-weak mixing decomposition states that for any unitary operator $U$ acting on a Hilbert space $\mathcal{H}$, we have a decomposition $\mathcal{H} = \mathcal{H}_c\oplus\mathcal{H}_w$ with

\begin{alignat*}{2}
    \mathcal{H}_c& = \{\xi \in \mathcal{H}\ :\ \xi\text{ has pre-compact orbit under }U\} = \overline{\text{lin}}\left(\left\{\xi \in \mathcal{H}\ :\ U\xi = \lambda\xi\text{ for some }\lambda \in \mathbb{S}^1\right\}\right),\\
    \mathcal{H}_w& = \left\{\xi \in \mathcal{H}\ :\ \lim_{N\rightarrow\infty}\frac{1}{N}\sum_{n = 1}^N|\langle U^n\xi,\xi\rangle| = 0\right\} = \left\{\xi \in \mathcal{H}\ :\ \lim_{N\rightarrow\infty}\frac{1}{N}\sum_{n = 1}^N|\langle U^n\xi,\eta\rangle| = 0\ \forall\ \eta \in \mathcal{H}\right\}.
\end{alignat*}
Koopman and von Neumann \cite{koopman1932dynamical} were the first to show that a dynamical system is weakly mixing if and only if it has no nontrivial eigenfunctions.
A version of the compact-weak mixing decomposition for locally compact groups was found by Godement \cite{GodementDecomposition} (see also \cite[Appendix 10.D]{TempelmanBook}), and a version for amenable semigroups acting by isometries was found by Dye \cite{dye1965ergodic}.
In between Godement and Dye, the compact-weak mixing decomposition was extended to various semigroups of operators on a Banach space by Jacobs, de Leeuw, and Glicksberg \cite{JacobsJdLG,DeLeeuwGlicksberg}, and is now known as the Jacobs-de Leeuw-Glicksberg (JdLG) decomposition. 
Several applications of the JdLG-decomposition can be found in the literature; for example, it is used in ergodic Ramsey theory \cite[Theorem 2.3]{ERTAnUpdate}, in the first proof of the Erd\H{o}s sumset conjecture \cite{TheErdosSumsetPaper}, in the Fürstenberg-Zimmer structure theorem \cite{EdekoKreidlerHaase}, in a proof of the Wiener-Wintner theorem \cite[Proposition 2.2]{assani2003wiener}, and more generally in the study of weighted ergodic theorems \cite{OnModulatedErgodicTheorems,MickySohail}. 
Further uses include results related to the entangled ergodic theorem \cite{entangled}, 
the Perron-Frobenius theory for completely positive operators on $W^*$-algebras \cite{Btkai2011Decomposition}, 
and evolution equations \cite[Section~V.2]{EN}.

A modern treatment of the JdLG-decomposition is given in \cite[Chapter 16]{OTAoET}.
After numerous results, it has been realized that as soon as one deals with unitary, or more generally isometric operators, splitting subspaces can be described as ${\cal H}_c$, ${\cal H}_w$ above, and under reasonable assumptions one does not even need the Banach space $\cal H$ to be a Hilbert space. However, for even the simplest non-isometric semigroups this decomposition can fail, see Example \ref{ex-non-isometric}.

To go beyond, one needs to define splitting subspaces differently. Let $E$ be a Banach space, and let $\mathscr{T}$ be a semigroup of bounded operators on $E$. Let $\mathscr{S}$ be the weak closure of $\mathscr{T}$, and set
    \begin{alignat}{2}\label{Er-Eaws-definition}
        E_r& = \{\xi \in E\ :\ \forall\ u \in \mathscr{S}, \exists\ v \in \mathscr{S}\text{ with }vu\xi = \xi\}, \\
        E_{\text{aws}}& = \{\xi \in E\ :\ u\xi = 0\text{ for some }u \in \mathscr{S}\}.
        \notag
    \end{alignat}
If $E$ is a Hilbert space and $\mathscr{T}$ is generated by a single unitary operator, we get the same subspaces as before. But in general, these sets are not linear subspaces and $E_{\text{aws}}$ is not $\mathscr{T}$-invariant, see Examples \ref{E_aws-not-a-subspace} and \ref{E_r-not-a-subspace}.


Another point is that ergodic theorems require in general bounded orbits, though similar decompositions may appear in other contexts. A natural framework is therefore to assume that the operator semigroup is relatively weakly compact; that is, that the orbit of every vector is relatively compact in the weak topology.

Now, a semigroup $\mathscr{T}$ of operators on a Banach space $E$ is said to be \textbf{JdLG-decomposable} if its weak closure $\mathscr{S}$ is weakly compact and we have a decomposition into closed $\mathscr{S}$-invariant subspaces $E = E_r\oplus E_{\text{aws}}$, with $E_r$ and $E_{\text{aws}}$ defined in \eqref{Er-Eaws-definition}. $E_r$ is called the \textit{reversible part of $E$} and $E_{\text{aws}}$ is the \textit{almost weakly stable part of $E$}.


Exact conditions of decomposability are as follows:

\begin{theorem}[cf. {\cite[Theorem 4.11]{DeLeeuwGlicksberg}, \cite[Theorem 16.5]{OTAoET}} ]\label{th-decomposable}
    Let $E$ be a Banach space and let $\mathscr{S}$ be a weakly compact semigroup of operators on $E$. Then $\mathscr{S}$ is JdLG-decomposable if and only if one of the following equivalent conditions holds:
    \begin{itemize}
    \item[(i)] the space $C(\mathscr{S})$ of continuous (bounded) functions on $\mathscr{S}$ admits a two-sided invariant mean;
    \item[(ii)] $\mathscr{S}$ has a {\it unique} minimal idempotent.
    \end{itemize}
\end{theorem}
In \cite[Chapter 16.3]{OTAoET}, a semigroup $\mathscr{T} \subseteq \mathcal{L}(E)$ is called JdLG-admissible if it is relatively compact and $\mathscr{S}:=\operatorname{cl}_\sigma(\mathscr{T})$ contains a unique minimal idempotent; thus, JdLG-admissible is the same as JdLG-decomposable.

The properties above are of course hard to check. The question is therefore to find other, more tractable conditions which imply (i) or (ii) of Theorem \ref{th-decomposable}. Here is a list of main known cases of JdLG-admissible semigroups. 

\begin{theorem}[cf. {\cite[Chapter 16.3]{OTAoET}}]\label{KnownExamplesOfJdLG}
    Let $E$ be a Banach space and $\mathscr{T}$ a semigroup of bounded linear operators on $E$.
    Then $\mathscr{T}$ is JdLG-admissible in the following cases:
    \begin{enumerate}[(i)]
    \item $\mathscr{T}$ is abelian and relatively weakly compact.

    \item $E$ is a Hilbert space and $\mathscr{T}$ consists of contractions.

    \item $E = L^1(X,\mu)$ with $(X,\mu)$ a probability space, and $\mathscr{T}$ consists of Markov operators.

    \item $E$ and its dual space $E'$ are both strictly convex, and $\mathscr{T}$ is relatively weakly compact.

    \item $E = L^p(X,\mu)$ for some measure space $(X,\mu)$ and $1 < p < \infty$, and $\mathscr{T}$ consists of contractions.
\end{enumerate}
\end{theorem}

\begin{remark}
\begin{enumerate}
    \item  In many cases of interest, the Banach space $E$ is reflexive.
    Consequently, it is worth recalling that a set of bounded linear operators on a reflexive Banach space $E$ is relatively weakly compact if and only if it is uniformly bounded.
    \item A special case of (i) is provided by monothetic operator semigroups. Let $E=L^1(X,\mu)$, where $(X,\mu)$ is a finite measure space and let $\mathscr{T}$ be generated by a contraction $T$ whose linear modulus $|T|$ is mean ergodic (e.g. if $T$ is a Dunford–Schwartz operator). Then $\mathscr{T}$ is JdLG-admissible (cf. \cite[Corollary 1.3]{KornfeldLin2000}).
\end{enumerate}

\end{remark}

Let now $S$ be a semigroup which is not a priori realized as a semigroup of operators.
We call a representation $\pi$ of $S$ on a Banach space $E$ {\it relatively weakly compact } if $\mathscr{T} := \{\pi_s: s\in S\}$ is relatively weakly compact (for example, $E$ reflexive and $\mathscr{T}$ uniformly bounded).
In this case, every function $s\mapsto f^\prime(\pi_s\xi)$, with $\xi \in E$ and $f^\prime \in E^\prime$, belongs to the space $W(S)$ of weakly almost periodic functions on $S$ (cf. \cite[Theorem 2.7]{DeLeeuwGlicksberg}). Moreover, $\pi$ is said to be \textbf{JdLG-admissible} if $\pi(S)$ is JdLG-admissible.
We have then:

\begin{theorem}[cf. {\cite[Theorem 6.2.18]{AnalysisOnSemigroups}}]\label{GeneralJdLGDecomposition}
    Let $S$ be a semigroup which admits a bi-invariant mean on $W(S)$, the space of weakly almost periodic functions on $S$.
    Let $E$ be a Banach space and let $\pi$ be a relatively weakly compact  representation of $S$ on $E$.
Then $\pi(S)$ is JdLG-admissible.
\end{theorem}

If $S$ is a topological group, or a bi-amenable measured semigroup, then it has a unique invariant mean on $W(S)$.
Consequently, the next result follows immediately from Theorem \ref{GeneralJdLGDecomposition}.
\begin{corollary}\label{NewAdditionsToJDLG}
    Let $S$ be either a topological group or a bi-amenable measured semigroup.
    If $E$ is a Banach space and $\pi$ is a relatively weakly compact representation of $S$ on $E$, then $\pi$ is JdLG-admissible.
    In particular, if $\mathscr{T}$ is a semigroup of bounded linear operators on $E$ with $c\ell_\sigma(\mathscr{T}) = c\ell_\sigma(\pi(S))$, then $\mathscr{T}$ is JdLG-admissible.
\end{corollary}


A powerful feature of the JdLG decomposition is the availability of multiple characterizations of the splitting subspaces. In the core of the paper, we recall several known descriptions and add a few new ones to them.

Suppose that $\mathscr{T} \subset {\cal L}(E)$ is a JdLG-admissible semigroup and $\mathscr{S}$ its weak closure.
On the reversible subspace $E_r$, the semigroup $\mathscr{S}$ acts as a compact group. This implies the existence of unitary finite-dimensional subrepresentations \cite[Theorem 16.31]{OTAoET}, and compactness of orbits \cite[Proposition 16.29]{OTAoET}.

The almost weakly stable subspace $E_{\text{aws}}$ consists of ``flight vectors'', a term given by Jacobs; every vector in $E_{\text{aws}}$ ``tends to zero in mean''. In the case of a single unitary operator acting on a Hilbert space, a precise meaning of this is given in the definition of $ \mathcal{H}_w$. In general, this convergence can be expressed in terms of the invariant mean $M$ on $C(\mathscr{S})$, which does exist by Theorem \ref{th-decomposable}.

We prove in particular the following characterization, where (i) is Theorem \ref{(Weak)EquivalenceForReversiblePart} below, and (ii) is a part of Lemma \ref{MeanCharacerizationOfAWS}:

\begin{theorem}\label{NewCharacterizationsOfRAndAWS}
    Let $E$ be a Banach space, $\mathscr{T} \subseteq \mathcal{L}(E)$ a JdLG-admissible semigroup of operators, $\mathscr{S}$ its weak closure, and $E = E_r\oplus E_{\text{aws}}$ the JdLG-decomposition.
    \begin{enumerate}[(i)]
        \item The restriction of $\mathscr{T}$ to $E_r$ is weakly equivalent to a unitary representation $U$ of $\mathscr{T}$ on a Hilbert space $\mathcal{H}$.
        Furthermore, the representation $U$ decomposes into a direct sum of finite dimensional representations.

        \item  Let $M$ be a bi-invariant mean on $C(\mathscr{S})$. 
        We have
        \begin{equation}
            E_{\text{aws}} = \left\{\xi \in E\ :\ \forall\ f^\prime \in E^\prime, M\left(T\mapsto \left|f^\prime(T\xi)\right|\right) = 0\right\}.
        \end{equation}
    \end{enumerate}
\end{theorem}



This is an analogue of the decomposition of weakly almost periodic functions (cf. \cite[Theorem 5.7, Corollary 5.9]{DeLeeuwGlicksberg} and \cite[Theorem 4.3.13]{AnalysisOnSemigroups}):
It is known that a semigroup $S$ admits a bi-invariant mean $m$ on $W(S)$ if and only if we have $W(S) = \mathcal{SAP}(S)\oplus W_0(S)$, where $\mathcal{SAP}(S)$ is the space of strongly almost periodic functions on $S$, and $W_0(S) = \{f \in W(S)\ :\ m(|f|) = 0\}$.
The space $\mathcal{SAP}(S)$ is also the closure of the space of coefficients of continuous finite-dimensional unitary representations of $S$, and $W_0(S)$ consists of those functions $f \in W(S)$ that have $0$ in the weak closure of the orbit of their right translates.
We see that this is a special case of a JdLG decomposition. 




When $S = \mathbb{N}$ in Theorem \ref{GeneralJdLGDecomposition}, multiple alternative characterizations of $E_{\text{aws}}$ resembling the characterizations of $\mathcal{H}_w$ are stated in \cite[Theorem 16.34]{OTAoET}.
In Theorem \ref{EquivalentCharacterizationsOfAWSWhenAmenable} we extend these alternative characterizations to various classes of amenable semigroups.
For now, we only state a particularly nice result following from Corollary \ref{NewAdditionsToJDLG}, Corollary \ref{ReversiblePartIsRelativelyCompactOrbits}, and Theorem \ref{EquivalentCharacterizationsOfAWSWhenAmenable}.


\begin{theorem}\label{intro-chars-by-foelner}
    Let $G$ be a locally compact amenable group, $\lambda$ a left-Haar measure, $\mathcal{F}=(F_i)_{i \in I}$ a left-F\o{}lner net, $E$ a Banach space, and $\pi$ a relatively weakly compact representation of $G$ on $E$.
    Then $\pi$ is JdLG-admissible, and we further have that
    \begin{alignat*}{2}
        E_{r}&=\{\xi\in E: \pi(G)\xi \text{ is relatively compact in }E\}\text{, and}\\
        E_{\text{aws}}&=\left\{\xi\in E:\lim_i\frac{1}{\lambda(F_i)}\int_{F_i}\left|\left\langle \pi_g \xi,x^\prime \right\rangle\right|d\lambda(g)=0 \text{ for all }x^\prime\in E^\prime\right\}\\
        &=\left\{\xi\in E:\lim_i\sup_{\|x^\prime\| = 1}\frac{1}{\lambda(F_i)}\int_{F_i}\left|\left\langle \pi_g \xi,x^\prime \right\rangle\right|d\lambda(g) =0 \right\}.
    \end{alignat*}
\end{theorem}

Results relating the compact-weak mixing decomposition to the algebra of ultrafilters on $\mathbb{N}$ have been investigated in \cite{bergelson2007central,BMIPUltrafiltersAndSzemerediForGPolynomials,bergelson2021iterated,mccutcheon2014d}, and in Section \ref{UltrafilterSection} we show that many of these results extend to JdLG-admissible semigroups of operators on a Banach space.

The structure of this paper is as follows. 
In Section 2, we collect known facts about topological semigroups, weak almost periodicity, JdLG-admissible semigroups, amenability, representations of semigroups, and ultrafilters.
In Section 3 we study the reversible part of the JdLG-decomposition.
In Section 4 we study the almost weakly stable part of the JdLG-decomposition.
In Section \ref{UltrafilterSection} we study the relationship between the JdLG-decomposition and the algebra of the Stone-\v{C}ech compactification.

\textbf{Acknowledgements:} The second author acknowledges being supported by grant 2019/34/E/S$\allowbreak$T1$\allowbreak$/00082 for the project “Set theoretic methods in dynamics and number theory,” NCN (The
National Science Centre of Poland), and the first author was supported by the same grant for two research visits.
The second author would like to thank Uta Freiberg for funding a research visit to TU Chemnitz.
We would also like to thank Rigoberto Zelada for helpful discussions about Theorem \ref{ProjectionsViaUltrafilters}.
Finally, we would like to thank Chat gpt 5.6 sol for proof reading this paper, and helping improve some of the results of a previous edition of this work.

\section{Preliminaries}
\subsection{Topological semigroups and weak almost periodicity}\label{WAPSubsection}
All of the facts that we collect in this section about weakly almost periodic functions come from the classical paper of de Leeuw and Glicksberg \cite{DeLeeuwGlicksberg}.
For further discussions on this topic the reader is referred to \cite{AnalysisOnSemigroups,burckel1970weakly}.

A semigroup $(S,\cdot)$ is a set $S$ with an associative binary operation $\cdot:S\times S\rightarrow S$.
If $S$ possesses an identity element $e$, then $S$ is called a \textbf{monoid}.
If $S$ is also a Hausdorff topological space and the map $\cdot$ is continuous, then $S$ is a \textbf{topological semigroup}\footnote{We warn the reader that in \cite{DeLeeuwGlicksberg} and \cite{burckel1970weakly}, the term topological semigroup refers to what we call a semitopological semigroup.}.
If $S$ is a Hausdorff topological space for which $\cdot$ is separately continuous, i.e., for each $s \in S$ the maps $L_s(t) := s\cdot t$ and $R_s(t) := t\cdot s$ are continuous, then $S$ is a \textbf{semitopological semigroup}.
For the rest of this paper all semigroups are assumed to be semitopological unless it is explicitly stated otherwise.
By abuse of notation, we let $R$ act on $C_b(S)$ by $R_s(f)(t)=f(R_s(t))=f(ts)$, and similarly for $L$. If $\lambda$ is a Borel measure on $S$ such that 
\[\lambda(N)=0 \Longrightarrow \lambda(R_s^{-1}(N))=0\]
for every $s\in S$, then the same formula defines $R_s$ on $L^\infty(S,\lambda)$, and analogously for $L$. 
A function $f \in C_b(S)$ is \textbf{weakly almost periodic} if $\{R_sf\ :\ s \in S\}$ is relatively weakly compact.
We let $W(S)$ denote the space of weakly almost periodic functions on $S$.
If $S$ is a compact semigroup, then $W(S) = C(S)$.
For $X \in \{C_b(S),W(S),L^\infty(S,\lambda)\}$, a \textbf{left-invariant mean} on $X$ is a positive linear functional $m:X\rightarrow\mathbb{C}$ of norm $1$ that satisfies $m(L_sf) = m(f)$ for all $s \in S$ and all $f \in X$, and a right-invariant mean is defined similarly.
The semigroup $S$ admits a \textbf{bi-invariant mean} on $X$ if there exists a left-invariant mean $m$ on $X$ that is also a right-invariant mean.


If $W(S)$ admits a bi-invariant mean, then it is unique.
It is also known that $W(S)$ admits a bi-invariant mean if and only if it admits a left-invariant mean and a (potentially different) right-invariant mean.
Now let $S$ be a semigroup with a bi-invariant mean $m$ on $W(S)$.
If $\mathscr{S}$ is a compact semigroup and $\pi:S\rightarrow\mathscr{S}$ is a continuous homomorphism with dense image, then the map $\widetilde{\pi}:C(\mathscr{S})\rightarrow W(S)$ given by $(\widetilde{\pi}f)(s)=f(\pi_s)$ is well defined, and $M(f):=m(\widetilde{\pi}f)$ is the bi-invariant mean on $W(\mathscr{S})=C(\mathscr{S})$.
There exists a compact semigroup $S^w$, \textbf{the weakly almost periodic compactification of $S$}, and a canonical continuous homomorphism $i:S\rightarrow S^w$ with dense image, such that for any $f \in W(S)$ there is a unique $f^w \in W(S^w)$ satisfying $\tilde{i}f^w(s) = f(s)$.
By abuse of notation, we let $i:W(S)\rightarrow W(S^w)$ denote the map $f\mapsto f^w$, and we observe that $i$ is an algebra isomorphism.
It is known that $W(S^w)$ admits a bi-invariant mean $M$ if and only if $K(S^w)$, the smallest two-sided ideal of $S^w$, is a compact topological group, in which case we have $M(f) = \int_{K(S^w)}fd\mu$, where $\mu$ is the normalized Haar measure on $K(S^w)$.

Now suppose that $W(S^w)$ admits a bi-invariant mean $M$ and a left-invariant mean $M_\ell$.
For $f \in W(S^w)$ and $g \in K(S^w)$, we have $M_\ell(f) = M_\ell(L_gf)$, hence $M_\ell(f)$ is determined by the values of $f$ on $K(S^w)$. Tietze's Extension Theorem tells us that any $f \in C(K(S^w))$ extends to a member of $W(S^w)$, so the left-invariant means on $W(S^w)$ are in a one-to-one correspondence with the left-invariant means on $C(K(S^w))$.
The uniqueness of the Haar measure tells us that there is exactly one left-invariant mean on $C(K(S^w))$, hence $M_\ell = M$.
Since the means on $W(S)$ are in a one-to-one correspondence with the means on $W(S^w)$, the analogous result holds for $W(S)$.

A group $(G,\cdot)$ is a \textbf{topological group} if it is a topological semigroup and the inversion map is a continuous map.
All groups that we consider in this paper will be topological groups.
De Leeuw and Glicksberg \cite{DeLeeuwGlicksberg} investigated which topological groups $G$ admit a bi-invariant mean on $W(G)$, and $6$ years later, Ryll-Nardzewski \cite{RNFixedPoint} proved a fixed point theorem that allows one to deduce that all topological groups $G$ admit a bi-invariant mean on $W(G)$.

\subsection{JdLG-admissible semigroups}\label{JdLG-adimissableSubsection}
For a Banach space $E$ we let $\mathcal{L}(E)$ denote the space of bounded linear operators on $E$.
$\mathcal{L}(E)$ is naturally a semitopological semigroup when endowed with the operation of composition as well as the weak operator topology.
In \cite[Chapter 16.3]{OTAoET}, a semigroup $\mathscr{T} \subseteq \mathcal{L}(E)$ is called JdLG-admissible if it is relatively compact and $\mathscr{S}:=\operatorname{cl}_\sigma(\mathscr{T})$ contains a unique minimal idempotent $Q$.
Let us show in detail that this is equivalent to being JdLG-decomposable as defined in the introduction.
If $\mathscr{T}$ is a relatively weakly compact semigroup, then \cite[Theorem 16.5]{OTAoET} says that the following are equivalent: 

\begin{enumerate}
    \item $\mathscr{T}$ is JdLG-admissible.

    \item $K(\mathscr{S})$, the smallest two-sided ideal of $\mathscr{S}$, is a group.

    \item $\mathscr{S}$ has a unique minimal right ideal and a unique minimal left ideal.
\end{enumerate}
\cite[Theorem 16.3]{OTAoET} tells us that every minimal right ideal of $\mathscr{S}$ is closed, and \cite[Lemma 16.4]{OTAoET} tells us that $K(\mathscr{S})$ is the union of all minimal right ideals of $\mathscr{S}$.
It follows that when $\mathscr{T}$ is JdLG-admissible, $K(\mathscr{S})$ is a closed (hence compact) subsemigroup of $\mathscr{S}$.
\cite[Theorem 2.1]{DeLeeuwGlicksberg} says that a compact semitopological semigroup that is algebraically a group must be a topological group.
It follows that $\mathscr{T}$ is JdLG-admissible if and only if $K(\mathscr{S})$ is a compact topological group. 
\cite[Theorem 4.11]{DeLeeuwGlicksberg} states that $K(\mathscr{S})$ is a compact topological group if and only if $\mathscr{T}$ is JdLG-decomposable.
Putting everything together, we see that $\mathscr{T}$ is JdLG-decomposable if and only if it is JdLG-admissible.

\begin{example}\label{E_aws-not-a-subspace}
On $E=\C^2$, let $a$ and $b$ be the operators acting as $(x,y)\mapsto (x,x)$ and $(x,y)\mapsto (y,y)$ respectively. They generate the (weakly compact) semigroup $S=\{I,a,b\}$ with multiplication $ba=a^2=a$ and $ab=b^2=b$. The only nontrivial $S$-invariant subspace is $\C(1,1)$, and the set $E_{\text{aws}}$ in \eqref{Er-Eaws-definition} is equal to $\{0\}\times \C \cup \C\times \{0\}$, which is not even a linear subspace and is not $S$-invariant.
\end{example}

\begin{example}\label{E_r-not-a-subspace}
It is easy to check that $E_r$ in \eqref{Er-Eaws-definition} is always $\mathscr S$-invariant. But it might not be a linear subspace.
Dualizing the previous example, we get a semigroup $S^* = \{ I, a^*, b^*\}$ acting as $a^*(\xi,\eta) = (\xi+\eta,0)$ and $b^*(\xi,\eta) = (0,\xi+\eta)$. Any reversible vector must lie in the range of $a^*$ or $b^*$, so it is easy to see that $E_r = \{0\}\times \C \cup \C\times \{0\}$.
\end{example}

We conclude this subsection with several structural properties of
JdLG-admissible representations and their associated decompositions.

\begin{theorem}
    Let $S$ be a semigroup, $E$ a Banach space, and $\pi$ a JdLG-admissible representation of $S$ on $E$.
    Let $E = E_r\oplus E_{\text{aws}}$ be the JdLG-decomposition, and let $P_\pi:E\rightarrow E_r$ denote the projection whose kernel is $E_{\text{aws}}$.
    \begin{enumerate}[(i)]      
        \item If $Y$ is a closed $\pi$-invariant subspace of $E$, then the restriction $\pi|_Y$ of $\pi$ to $Y$ is again JdLG-admissible.
        Furthermore, we have $Y_r = Y\cap E_r$, $Y_{\text{aws}} = Y\cap E_{\text{aws}}$, and $P_{\pi|_Y} = (P_\pi)|_Y$.
        
        \item If $\tau$ is a representation of $S$ on the Banach space $F$, and $T:E\rightarrow F$ is a surjective bounded linear intertwining map for $\pi$ and $\tau$, then $\tau$ is JdLG-admissible.
        Furthermore, we have $F_r = TE_r$, $F_{\text{aws}} = TE_{\text{aws}}$, and $P_\tau T = TP_\pi$, with the obvious notation for $P_\tau$.

        \item If $Y$ is a closed $\pi$-invariant subspace of $E$, then the quotient representation $\pi^q$ of $S$ on the quotient space $E/Y$ is again JdLG-admissible.
        Furthermore, we have $(E/Y)_r = (E_r+Y)/Y$, $(E/Y)_{\text{aws}} = (E_{\text{aws}}+Y)/Y$, and $P_{\pi^q}$ is given by $P_{\pi^q}(\xi+Y) = P_\pi\xi+Y$, $\xi\in E$.

        \item If the adjoint representation $\pi'$ of the opposite semigroup $S^{\text{op}}$ on $E'$ is relatively weakly compact,\footnote{This condition is automatically fulfilled if $E$ is reflexive, a case considered already by Jacobs \cite{JacobsJdLG}.} then $(E')_r = E_{\text{aws}}^\perp$, $(E')_{\text{aws}} = E_r^\perp$, and $P_{\pi'} = P_\pi'$.

        \item The semigroup $\mathscr{S} := \text{cl}_\sigma(\pi(S))$ is a subdirect product of a group $G$ and a semigroup $R$ whose minimal two-sided ideal is a singleton.
    \end{enumerate}
\end{theorem}

\begin{proof}
    Part (i) is immediate, so we begin with the proof of (ii).
    Since $T$ is a surjective bounded linear intertwining map, the map

    \begin{equation}
        \Phi:\text{cl}_\sigma(\pi(S))\rightarrow\text{cl}_\sigma(\tau(S))
    \end{equation}
    given by $\Phi(u)(Tx) = Tux$, is a well-defined continuous surjective homomorphism.
    It follows that $\tau(S)$ is relatively weakly compact.
    Furthermore, if $m$ is a two-sided invariant mean on $C(\text{cl}_\sigma(\pi(S)))$, then $\Phi_*m$ is a two-sided invariant mean on $C(\text{cl}_\sigma(\tau(S)))$, so $\tau$ is indeed JdLG-admissible.
    We now recall that $P_\pi$ is the unique idempotent in $K(\text{cl}_\sigma(\pi(S)))$ and $P_\tau$ is the unique idempotent in $K(\text{cl}_\sigma(\tau(S)))$.
    Since $\Phi$ is a continuous surjective homomorphism, we see that $\Phi(P_\pi) = P_\tau$, i.e., that $P_\tau T = TP_\pi$.
    We now see that $F_r = P_\tau F = P_\tau TE = TP_\pi E = TE_r$, and that $F_{\text{aws}} = \text{ker}(P_\tau) = TT^{-1}\text{ker}(P_\tau) = T\text{ker}(P_\tau T) = T\text{ker}(TP_\pi) = T(\text{ker}(P_\pi)+\text{ker}(T)\cap\text{Im}(P_\pi)) = TE_{\text{aws}}$.

    We now see that (iii) follows from (ii) by letting $T:E\rightarrow E/Y$ be the quotient map.

    We now proceed to prove (iv).
    We see that adjunction is a homeomorphic anti-isomorphism from $\text{cl}_\sigma(\pi(S))$ to $\text{cl}_{\sigma^*}(\pi'(S^{\text{op}}))$, where $\sigma^*$ denotes the weak$^*$ operator topology.
    Since $\text{cl}_\sigma(\pi'(S^{\text{op}}))$ is compact by assumption we have that $\text{cl}_\sigma(\pi'(S^{\text{op}})) = \text{cl}_{\sigma^*}(\pi'(S^{\text{op}}))$, and that the weak operator and weak$^*$ operator topologies agree on this set.
    Consequently, adjunction maps the unique minimal idempotent of
$\operatorname{cl}_\sigma(\pi(S))$ to the unique minimal idempotent of
$\operatorname{cl}_\sigma(\pi'(S^{\mathrm{op}}))$, and hence
$P_{\pi'}=P_\pi'$. We now see that $(E')_r = \text{Im}(P_{\pi'}) = \text{ker}(P_\pi)^\perp = E_{\text{aws}}^\perp$ and $(E')_{\text{aws}} = \text{ker}(P_{\pi'}) = \text{Im}(P_\pi)^\perp = E_r^\perp$.

    We now proceed to prove (v).
    Let $\mathscr{S} = \text{cl}_\sigma(\pi(S))$, $G = \{u|_{E_r}\ |\ u \in \mathscr{S}\}$, and $R = \{u|_{E_{\text{aws}}}\ |\ u \in \mathscr{S}\}$.
    We see that the map $u\mapsto(u|_{E_r},u|_{E_{\text{aws}}})$ is a homomorphism from $\mathscr{S}$ to $G\times R$ whose image projects surjectively onto both coordinates, so $\mathscr{S}$ is a subdirect product of $G$ and $R$.
    The fact that $G$ is a group is immediate from the definition of the reversible part.
    As for $R$, it contains the zero operator, so $\{0\}$ is a two-sided ideal of $R$.
Moreover, every nonempty two-sided ideal of $R$ contains $0$, and hence $K(R)=\{0\}$.
\end{proof}

\subsection{Amenable semigroups and F\o lner nets}\label{AmenableSemigroupsSubsection}
The pair $(S,\lambda)$ is a \textbf{measured semigroup} if $S$ is a semitopological semigroup endowed with a Borel measure $\lambda$. 
The main cases of interest are when $S$ is a discrete semigroup and $\lambda$ is the counting measure, or when $S$ is a subsemigroup of a locally compact topological group and $\lambda$ is the restriction of a Haar measure.
The pair $(S,\lambda)$ is \textbf{left-amenable} if $\lambda(N)=0\Longrightarrow \lambda(L_s^{-1}N)=0$ for every $s\in S$ and every Borel set $N\subseteq S$, and there exists a left-invariant mean $m$ on $L^\infty(S,\lambda)$.
The pair $(S,\lambda)$ is \textbf{bi-amenable} if $\lambda(N)=0\Longrightarrow \lambda(L_s^{-1}N)=\lambda(R_s^{-1}N)=0$ for every $s\in S$ and every Borel set $N\subseteq S$, and there exists a bi-invariant mean on $L^\infty(S,\lambda)$.
For the rest of this section let $S$ denote a locally compact semitopological semigroup and let $\lambda$ be a Radon measure on $S$.
A \textbf{left-F\o lner net} $\mathcal{F} = (F_i)_{i \in I}$ is a net of compact subsets of $S$ of positive measure satisfying 

\begin{equation}
    \lim_i\left\|(L_s)_*\frac{\lambda|_{F_i}}{\lambda(F_i)}-\frac{\lambda|_{F_i}}{\lambda(F_i)}\right\|_{\text{TV}} = 0\text{ for all }s \in S.
\end{equation}
$\mathcal{F}$ is a \textbf{two-sided F\o lner net} if it is simultaneously a left-F\o lner net and a right-F\o lner net. 
In the discrete case with counting measure, for every finite nonempty 
$F\subseteq S$, putting
$\mu_F:=\frac{1}{|F|}\sum_{t\in F}\delta_t,$
we have
$\|(L_s)_*\mu_F-\mu_F\|_{\mathrm{TV}}
=2\frac{|F\setminus sF|}{|F|}.$
Thus, in the discrete setting, our definition of a left-F\o lner net
agrees with the usual strong F\o lner condition (SFC). If $S$ is left
cancellative, this condition can equivalently be expressed in terms of
$|sF\triangle F|/|F|$.

Suppose, in addition, that the corresponding translations preserve
$\lambda$-null sets under preimages. More precisely, assume that
\[
\lambda(N)=0\Longrightarrow \lambda(L_s^{-1}(N))=0
\]
in the left case,
\[
\lambda(N)=0\Longrightarrow \lambda(R_s^{-1}(N))=0
\]
in the right case, and that both conditions hold in the two-sided case,
for every $s\in S$ and every Borel set $N\subseteq S$.
If $S$ possesses a left/right/two-sided F\o lner net
$\mathcal{F}=(F_i)_{i\in I}$, then $(S,\lambda)$ is
left/right/bi-amenable, respectively. To see this, we define a net of functionals $(\mu_i)_{i \in I}$ on $L^\infty(S,\lambda)$ by $\mu_i(f) = \frac{1}{\lambda(F_i)}\int_{F_i}f(s)d\lambda(s)$ for all $f \in L^\infty(S,\lambda)$.
Since the unit ball of $L^\infty(S,\lambda)^\prime$ is weak$^*$ compact, we let $\mu$ be any cluster point of $(\mu_i)_{i\in I}$ and we see that $\mu$ is a left/right/two-sided invariant mean on $L^\infty(S,\lambda)$.
We say that $\mu$ is a mean generated by $\mathcal{F}$ and we denote the set of such means by $\mathcal{M}(\mathcal{F})$.

There are examples of discrete semigroups that are amenable but do not possess a F\o lner net.
However, if $S$ is a discrete cancellative semigroup, then $S$ is amenable if and only if it possesses a F\o lner net.
It is also well known that if $S$ is a discrete abelian semigroup, then $S$ possesses a two-sided F\o lner net.
We refer the reader to \cite[Chapter 4.22]{PatersonAmenability} for further discussion on F\o lner conditions and amenability for discrete semigroups.

Now suppose that $S$ is a left-amenable semigroup.
We define the upper (left) and lower (left) Banach densities on the collection of measurable subsets of $S$ by
\begin{alignat*}{2}
    d^*(A)& = \sup\{m(A) | m\text{ is a left-invariant mean on }S\},\text{ and}\\
    d_*(A)& = \inf\{m(A) | m\text{ is a left-invariant mean on }S\} = 1-d^*(A^c).
\end{alignat*}
If we further assume that $S$ possesses a left-F\o lner net $\mathcal{F} = (F_i)_{i \in I}$, then we define the upper and lower densities along $\mathcal{F}$ by
\[\overline{d}_{\mathcal{F}}(A):=\limsup_i\frac{\lambda(A\cap F_i)}{\lambda(F_i)}\text{ and }\underline{d}_{\mathcal{F}}(A):=\liminf_i\frac{\lambda(A\cap F_i)}{\lambda(F_i)},\]
respectively.
If $\overline{d}_{\mathcal{F}}(A) = \underline{d}_{\mathcal{F}}(A)$, then we denote the common value by $d_{\mathcal{F}}(A)$.
We define the upper (left) and lower (left) F\o lner densities on $S$ by
\begin{alignat*}{2}
    \overline{d}(A)& = \sup\left\{\overline{d}_\mathcal{F}(A) | \mathcal{F}\text{ is a (left) F\o lner net in }S\right\},\text{ and}\\
    \underline{d}(A)& = \inf\left\{\underline{d}_\mathcal{F}(A) | \mathcal{F}\text{ is a (left) F\o lner net in }S\right\} = 1-\overline{d}(A^c).
\end{alignat*}
If $S$ is a discrete semigroup and $\lambda$ is the counting measure, then \cite[Theorem 3.15]{glasscock2025folner} tells us that $\overline{d} = d^*$ and $\underline{d} = d_*$.
More generally, we see that $\overline{d} \le d^*$ (hence $\underline{d} \ge d_*$) since any F\o lner net generates many means, but it is not clear if $d^* \le \overline{d}$.

Suppose that we are given a measurable function $f:S\rightarrow \mathbb{C}$ and $D \in \{d_*,\underline{d},d_{\mathcal{F}}\}$ that is well-defined with respect to $S$.
We write
\[D-\lim_{s} f(s) = \alpha\]
if for any open neighborhood $U$ of $\alpha$, 
\[D(\{s \in S\ :\ f(s)\in U\})=1.\]
Given a topological space $Y$, a measurable function $f:S\rightarrow Y$, and a set $A \subseteq S$, we write 
\[\lim_{\underset{s \in A}{s\to \infty}}f(s) = \alpha\]
if for any open neighborhood $U$ of $\alpha$ the set $\{s \in A\ :\ f(s) \notin U\}$ is contained in a compact set.

\subsection{Representations of semigroups}\label{RepresentationsOfSemigroupsSubsection}
Let $E$ be a Banach space, $\mathcal{L}(E)$ the space of bounded linear operators on $E$, and $S$ a semigroup.
A \textbf{representation} $\pi$ of $S$ on $E$ is a weakly continuous\footnote{We warn the reader that many other texts require that $\pi$ be strongly continuous. Nonetheless, for our purposes weak continuity is enough.} semigroup homomorphism $\pi:S\rightarrow\mathcal{L}(E)$.
If $S$ is a monoid with identity $e$, then we will assume that $\pi_e$ is the identity map.
We note that if $\mathcal{H}$ is a Hilbert space and $U$ is a unitary representation of $S$ on $\mathcal{H}$, then $U$ is strongly continuous.
The representation $\pi$ is \textbf{uniformly bounded} if $\sup_{s \in S}\|\pi_s\| < \infty$, and it is \textbf{relatively weakly compact} if $\{\pi_s\ :\ s \in S\}$ is relatively compact in the weak topology of $\mathcal{L}(E)$.
A fact that we will frequently use is that if $\pi$ is a relatively weakly compact representation of $S$ on $E$, then for any $\xi \in E$ and any $f^\prime \in E^\prime$, the function $\phi(s) = f^\prime(\pi_s\xi)$ belongs to $W(S)$ (cf. \cite[Theorem 2.7]{DeLeeuwGlicksberg}). 
Since $W(S)$ is a $C^*$-algebra, it follows that $|\phi|^p\in W(S)$ as well for every $p > 0$.
While every relatively weakly compact set of operators in $\mathcal{L}(E)$ is uniformly bounded, the converse need not be true.

Let $\mathcal{I}(E)$ be the set of invertible bounded linear operators on $E$.
If $H$ is a uniformly bounded subgroup of $\mathcal{I}(E)$, then $H$ is a topological group when given the operator norm topology.
Recently, Banakh \cite{banakh2022automatic} showed that a Haar-measurable homomorphism\footnote{Given two topological groups $X$ and $Y$ with $X$ locally compact, a map $h:X\rightarrow Y$ is Haar-measurable if for every open set $U \subseteq Y$, the set $h^{-1}(U)$ is in the completion of the Borel $\sigma$-algebra of $X$ with respect to a Haar measure of $X$.} from a locally compact group to a topological group is in fact a continuous homomorphism.
In particular, if $G$ is a locally compact group and $\pi:G\rightarrow\mathcal{I}(E)$ is a relatively weakly compact Haar-measurable homomorphism, then $\pi$ is strongly (hence weakly) continuous.

A \textbf{matrix coefficient} $\phi$ of the representation $\pi$ is a function of the form $\phi(s) = f^\prime(\pi_s\xi)$ for some $\xi \in E$ and some $f^\prime \in E^\prime$.
If $\pi_1$ and $\pi_2$ are representations of $S$ on $E_1$ and $E_2$ respectively, then $\pi_1$ is \textbf{weakly contained} in $\pi_2$ if every matrix coefficient of $\pi_1$ can be approximated arbitrarily well by linear combinations of matrix coefficients of $\pi_2$.
To be more precise, $\pi_1$ is weakly contained in $\pi_2$ if for any $\xi \in E_1$, any $f^\prime \in E_1^\prime$, and any $\epsilon > 0$, there exist $\eta_1,\cdots,\eta_k \in E_2$ and $h^\prime_1,\cdots,h^\prime_k \in E^\prime_2$ for which

\begin{equation}
    \sup_{s \in S}\left|f^\prime(\pi_{1,s}\xi)-\sum_{j = 1}^kh^\prime_j(\pi_{2,s}\eta_j)\right| < \epsilon.
\end{equation}
The representations $\pi_1$ and $\pi_2$ are \textbf{weakly equivalent} if $\pi_1$ is weakly contained in $\pi_2$ and $\pi_2$ is weakly contained in $\pi_1$.
The representations $\pi_1$ and $\pi_2$ are \textbf{equivalent} if there exists a bounded linear invertible operator $T:E_1\rightarrow E_2$ for which $\pi_{2,s}T = T\pi_{1,s}$ for all $s \in S$.

The following classical observation will be useful for us later on.
If $\pi$ is a uniformly bounded representation of a monoid $S$ on $E$, then we can endow $E$ with a new norm $\|\cdot\|_\pi$ given by $\|\xi\|_\pi = \sup_{s \in S}\|\pi_s\xi\|$.
We observe that $\pi$ is a representation of $S$ by contractions on the Banach space $(E,\|\cdot\|_\pi)$.
If $S$ is a monoid, then the norm $\|\cdot\|_\pi$ is equivalent to $\|\cdot\|$.
It is worth recalling that replacing the norm of a Hilbert space with an equivalent norm may produce a Banach space that is not a Hilbert space.
In fact, from \cite{ehrenpreis1955uniformly} one can deduce that there exists a uniformly bounded representation $\pi$ of SL$(2,\mathbb{R})$ on a Hilbert space $\mathcal{H}$ such that $(\mathcal{H},\|\cdot\|_\pi)$ is not a Hilbert space.

\subsection{Ultrafilters and the algebra of the Stone-\v{C}ech compactification}\label{UltrafiltersSubsection}
If $S$ is a set, we denote by $\mathcal{P}(S)$ its power set.
\begin{definition}
    Given a set $S$, a collection of subsets $\mathcal{U}\subseteq \mathcal{P}(S)$ is called a \textbf{filter} if 
    \begin{enumerate}
        \item[(i)] $\emptyset\notin \mathcal{U}$ and $S\in \mathcal{U}$.
        \item[(ii)] If $A,B\in \mathcal{U}$ then $A\cap B\in \mathcal{U}$. 
        \item[(iii)] If $A\in \mathcal{U}$ and $A\subseteq B$ then $B\in \mathcal{U}$.
    \end{enumerate}
    The collection $\mathcal{U}$ is called an \textbf{ultrafilter} if it is a filter and for any $A \subseteq S$ we have either $A \in \mathcal{U}$ or $A^c \in \mathcal{U}$.
\end{definition}
We let $\beta S$ denote the space of ultrafilters on $S$, which is also known as the Stone-\v{C}ech compactification of $S$ when $S$ is endowed with the discrete topology.
We embed $S$ in $\beta S$ as a dense subset with the identification $s \mapsto e(s) := \{A \subseteq S\ |\ s \in A\}$, and we let $S^* := \beta S\setminus S$ denote the set of nonprincipal ultrafilters.
Given a topological space $X$, a function $f:S\rightarrow X$, and a point $x \in X$, we write $p-\lim_sf(s) = x$ if for any open neighborhood $U$ of $x$ we have $\{s \in S | f(s) \in U\} \in p$.
A useful property of ultrafilters is that for any compact Hausdorff topological space $X$, any $p \in \beta S$, and any function $f:S\rightarrow X$, there exists a unique $x \in X$ for which $p-\lim_sf(s) = x$.
We will make use of this fact for relatively weakly compact representations $\pi$ of a semigroup $S$ on a Banach space $E$, as we can then define for each $p \in \beta S$ the operator $\displaystyle\pi_p := p-\lim_s\pi_s$.

When $(S,\cdot)$ is a discrete semigroup, there is a unique extension of the operation $\cdot$ from $e(S)$ to all of $\beta S$ that makes $\beta S$ a compact \textit{right topological semigroup}, i.e., each of the maps $R_p$ for $p \in \beta S$ is continuous. 
One can just as well extend the semigroup operation from $S$ to $\beta S$ so as to make each of the maps $L_p$ continuous instead of $R_p$, but in general the operation cannot be extended in such a way as to make $\beta S$ a semitopological semigroup.
If $\mathscr{S}$ is a compact semitopological semigroup and $i:S\rightarrow\mathscr{S}$ is a continuous homomorphism whose image is dense in $\mathscr{S}$, then there exists a unique continuous surjective homomorphism $h:\beta S\rightarrow \mathscr{S}$ satisfying $h\circ e = i$, and it is given by $h(p) = p-\lim_si(s)$ (cf. \cite[Theorem 4.8]{AlgebraInTheSCC}).
It is known that $\beta S$ possesses the smallest two-sided ideal that we denote by $K(\beta S)$, and that $K(\beta S)$ possesses idempotents.
It turns out that $\overline{K(\beta S)}$ is a two-sided ideal.

If $S$ is a left-amenable\footnote{We mention that we work with a left-amenable semigroup instead of a right-amenable semigroup because we made $\beta S$ a right topological semigroup.} semigroup, then we let $\mathcal{D}^*(S) = \{p \in \beta S | \forall\ A \in p, d^*(A) > 0\}$\footnote{In \cite{glasscock2025folner} the set $\mathcal{D}^*(S)$ is denoted by $\Delta^*(S)$.
We choose not to use this notation because we will be using $\Delta^*$ as in \cite{bergelson2021iterated,QuotientSetsAndDensityRecurrentSets} instead.}.
In \cite[Theorem 2.8]{glasscock2025folner} it is shown that $\mathcal{D}^*(S)$ is a closed two-sided ideal in $\beta S$, hence it contains $\overline{K(\beta S)}$.
An idempotent element of $\mathcal{D}^*(S)$ is called an \textbf{essential idempotent}.
Essential idempotents play an important role in ergodic Ramsey theory (cf. \cite{BMIPUltrafiltersAndSzemerediForGPolynomials,mccutcheon2014d,beiglbock2009solvability,bergelson2016polynomial,campbell2018d}).

When $S = G$ is a group, we let $i:G\rightarrow G$ denote the inversion map.
For $p \in \beta G$, we write $p^{-1} := i(p) = p-\lim_gi(g) = p-\lim_gg^{-1}$.
For $A \subseteq G$, we write $A^{-1} := \{a^{-1}\ |\ a \in A\}$, and we see that $A \in p$ if and only if $A^{-1} \in p^{-1}$.

For further discussion about the algebra of $\beta S$ when $S$ is a discrete semigroup, the reader is referred to \cite{AlgebraInTheSCC}.
We mention that the algebra of Stone-\v{C}ech compactification of a semitopological semigroup $S$ has been investigated in \cite{baker1976stone}, but the situation there is much more complicated, which is why we choose to restrict our attention to discrete semigroups when working with the Stone-\v{C}ech compactification.

\subsubsection{Notions of largeness}\label{sssec:NotionsOfLargeness}
In this section we review many combinatorial notions of largeness in discrete semigroups that can also be characterized using ultrafilters.
To this end, let us fix a discrete semigroup $S$.
We will use $A$ and $B$ to denote subsets of $S$, and $\mathscr{P}_f(S)$ to denote the collection of all finite nonempty subsets of $S$.

The set $A$ is a \textbf{IP-set} if there exists a sequence $(s_n)_{n = 1}^\infty$ in $S$ such that for every $F \in \mathscr{P}_f(\mathbb{N})$ we have $\prod_{f \in F}s_f \in A$, where the product is taken in increasing order.
The set $A$ is a \textbf{IP$^*$-set} if for any IP-set $B$ we have $A\cap B \neq \emptyset.$
It is known that $A$ is an IP-set if and only if there exists an idempotent $p \in \beta S$ with $A \in p$ (cf. \cite[Theorem 5.12]{AlgebraInTheSCC}), and that $A$ is an IP$^*$-set if and only if for every idempotent $p \in \beta S$ we have $A \in p$ (cf. \cite[Theorem 16.6]{AlgebraInTheSCC}).

The set $A$ is \textbf{syndetic} if there exists $H \in \mathscr{P}_f(S)$ for which $\bigcup_{h \in H}h^{-1}A = S$.
The set $A$ is \textbf{thick} if for every $F \in \mathscr{P}_f(S)$ there exists $x \in S$ such that $Fx \subseteq A$.
The set $A$ is \textbf{piecewise syndetic} if there exists $H \in \mathscr{P}_f(S)$ for which the family $\{y^{-1}\bigcup_{h \in H}h^{-1}A\ |\ y \in S\}$ has the finite intersection property, and this is equivalent to $A$ being the intersection of a syndetic set and a thick set.
The set $A$ is \textbf{thickly syndetic} if for any piecewise syndetic set $B$ we have $A\cap B\neq \emptyset$.
It is known that $A$ is piecewise syndetic if and only if there exists a $p \in \overline{K(\beta S)}$ with $A \in p$ (cf. \cite[Theorem 4.40]{AlgebraInTheSCC}), and that $A$ is thickly syndetic if and only if for every $p \in \overline{K(\beta S)}$ we have $A \in p$.

Now let us assume that $S = G$ is a group.
The set $A$ is a \textbf{difference set} if there exists an injective sequence $(g_n)_{n = 1}^\infty$ in $G$ for which $\{g_n^{-1}g_m\ |\ n < m\} \subseteq A$.
It can be checked that $A$ is a difference set if and only if there exists $p \in S^*$ for which $A \in p^{-1}p$.
More generally, given $\ell \in \mathbb{N}$, the set $A$ is a \textbf{$\Delta^\ell$-set} if there exist injective sequences $(s_{n,1})_{n = 1}^\infty,\dots,(s_{n,\ell})_{n = 1}^\infty$ for which $\{s_{n_1,1}^{-1}s_{m_1,1}\dots s_{n_\ell,\ell}^{-1}s_{m_\ell,\ell}\ |\ n_1 < m_1 < \dots < n_\ell < m_\ell\} \subseteq A$.
It is known that $A$ is a $\Delta^\ell$-set if and only if there exist $p_1,\dots,p_\ell \in S^*$ for which $A \in p_1^{-1}p_1\dots p_\ell^{-1}p_\ell$ (cf. \cite[Theorem 5.4]{QuotientSetsAndDensityRecurrentSets}).
The set $A$ is a \textbf{$\Delta^{\ell*}$-set} if for any $\Delta^\ell$ set $B$ we have $A\cap B \neq \emptyset$.
It is known that $A$ is a $\Delta^{\ell*}$-set if and only if for any $p_1,\cdots,p_\ell \in S^*$ we have $A \in p_1^{-1}p_1\cdots p_\ell^{-1}p_\ell$ (cf. \cite[Corollary 5.5]{QuotientSetsAndDensityRecurrentSets}).
We will see in the next subsection that the preceding discussion about $\Delta^\ell$-sets extends to the setting of cancellative left reversible semigroups.

\subsubsection{Cancellative left reversible semigroups}

A semigroup $S$ is called \textbf{cancellative} if, for all $a,b,c\in S$,
$ab=ac$ implies $b=c$, and $ba=ca$ implies $b=c$. The semigroup $S$ is \textbf{left reversible} if any two right ideals of $S$ intersect, i.e., if for any $r,s \in S$ we have $rS\cap sS \neq \emptyset$.
A group $G$ is a \textbf{group of right quotients of $S$} if $S \subseteq G$ and every $g \in G$ has the form $g = st^{-1}$ for some $s,t \in S$. 
It is a classical result of Ore (cf. \cite[Theorem 1.23]{AlgebraicTheoryOfSemigroupsVol1}) that if $S$ is a cancellative left reversible semigroup, then $S$ embeds in a group of its right quotients, and in this case $S$ will be a thick subset of $G$.
For example, $(\mathbb{Z},+)$ is the group of right quotients of $(\mathbb{N},+)$, and $(\mathbb{Q}^\times,\cdot)$ is the group of right quotients of $(\mathbb{Z}\setminus\{0\},\cdot)$.
It is well known (cf. \cite[Proposition 1.23]{AmenabilityByPaterson}) that every left-amenable semigroup is left reversible, so Example \ref{CounterExampleForAWSCharacterizations} contains an example of a noncommutative cancellative left reversible semigroup. 
We will also be in need of the following result.

\begin{lemma}\label{lem:ThicklySyndeticInReversibleSemigroups}
    Let $S$ be a left reversible semigroup.
    For every $s \in S$, the set $S\setminus sS$ is not piecewise syndetic, hence $sS$ is thickly syndetic.
\end{lemma}

\begin{proof}
Let \(A=S\setminus sS\), and let \(F\subseteq S\) be finite and nonempty. For each \(t\in F\), left reversibility gives an element \(r_t\in S\) such that
\[
tr_t\in sS.
\]
Again using left reversibility, the finitely many right ideals \(r_tS\), \(t\in F\), have nonempty intersection. Choose
\[
r\in\bigcap_{t\in F}r_tS.
\]
Then \(tr\in sS\) for every \(t\in F\). Since \(sS\) is a right ideal, it follows that
\[
trx\in sS
\]
for every \(t\in F\) and \(x\in S\). Therefore,
\[
rx\notin\bigcup_{t\in F}t^{-1}A
\]
for every \(x\in S\). Thus
\[
\bigcup_{t\in F}t^{-1}A
\]
is not thick. Since this holds for every finite nonempty \(F\subseteq S\), the set \(A=S\setminus sS\) is not piecewise syndetic.
\end{proof}

Now let us assume that $S$ is a discrete cancellative left reversible semigroup with group of right quotients $G$.
We see that the natural inclusion map $\iota:S\rightarrow G$ has a unique continuous extension to a natural inclusion map $\tilde{\iota}:\beta S\rightarrow \beta G$ given by $\tilde{\iota}(p) = \{A \subseteq G\ |\ A\cap S \in p\}$.
For $q \in \text{Im}(\tilde{\iota})$ we have $\tilde{\iota}^{-1}(q) = \{A\cap S\ |\ A \in q\}$.
Henceforth we identify $\beta S$ as a subset of $\beta G$ without explicitly mentioning the map $\iota$. For $p_1,\dots,p_\ell\in S^*$, the product
$q:=p_1^{-1}p_1\cdots p_\ell^{-1}p_\ell$
is understood as a product in $\beta G$.
We call $A\subseteq S$ a $\Delta^\ell$-set if there exist
$p_1,\dots,p_\ell\in S^*$ such that $q\in\beta S$ and $A\in q$.
We call $A\subseteq S$ a $\Delta^{\ell*}$-set if $A\in q$
for every $p_1,\dots,p_\ell\in S^*$ for which $q\in\beta S$.

Since $S$ is thick in $G$, $\beta S$ contains a left ideal of $\beta G$, hence it also contains a minimal left ideal of $\beta G$.
Since $K(\beta G)$ is the disjoint union of all minimal left ideals of $\beta G$, we see that $\beta S\cap K(\beta G) \neq \emptyset$, hence $K(\beta S) = \beta S\cap K(\beta G)$ and $\overline{K(\beta S)} = \beta S\cap\overline{K(\beta G)}$.
In light of Lemma \ref{lem:ThicklySyndeticInReversibleSemigroups}, we see that for any $p \in \overline{K(\beta S)}$ and $s \in S$, we have $sS \in p$, hence $S \in s^{-1}p$.
It follows that $s^{-1}p \in \beta S$, and $s^{-1}p \in \overline{K(\beta G)}$, hence $s^{-1}p \in \overline{K(\beta S)}$.
We now see that for any $q \in \beta S$ and $p \in \overline{K(\beta S)}$ we have $q^{-1}p = q-\lim_ss^{-1}p \in \overline{K(\beta S)}$, i.e., $\overline{K(\beta S)}$ is a left ideal in $\beta G$.

Now let us further assume that $S$ is left-amenable.
We see that any left-invariant mean $m$ on $\ell^\infty(S)$ naturally extends to a left-invariant mean on $\ell^\infty(G)$ given by $m(f) = m(f|_S)$.
It follows that $\mathcal{D}^*(S) \subseteq \mathcal{D}^*(G)$, and that $\mathcal{D}^*(S)$ is a left ideal in $\beta G$.
It follows that for $q \in \beta S$ and $p \in \mathcal{D}^*(S)$, we have $q^{-1}p \in \mathcal{D}^*(S)$. 

We warn the reader that for arbitrary $p,q \in \beta S$, we need not have $q^{-1}p \in \beta S$.
Nonetheless, it can still be checked that if $p \in \beta S$ is such that $p^{-1}p \in \beta S$, then for every $A \in p^{-1}p$, $A\cap S$ contains a difference set.
It is still the case that for any difference set $A \subseteq S \subseteq G$, there exists $p \in \beta S$ with $p^{-1}p \in \beta S$ and $A \in p^{-1}p$.
The situation for $\Delta^\ell$-sets is identical.

The following universal property of the group of right quotients will also be of use later on.

\begin{proposition}[{\cite[Theorem 2.2.4]{RandomPaperForOreLemma}}]\label{prop:ExtendingSemigroupHomomorphisms}
Let $S$ be a cancellative left reversible semigroup with group of right quotients $G$.
If $G_2$ is a group and $h\colon S\to G_2$ is a semigroup homomorphism, then $h$ extends uniquely to a group
homomorphism $\widetilde h\colon G\to G_2$, given by $\widetilde h(ab^{-1})=h(a)h(b)^{-1}$.
\end{proposition}

\section{The reversible part of the JdLG-decomposition}\label{ReversiblePartSection}
As in Section \ref{JdLG-adimissableSubsection}, we denote by $\mathscr{T}$ a relatively weakly compact operator semigroup acting on a Banach space $E$, and we denote by $\mathscr{S}$ the weak closure of $\mathscr{T}$. 
We have already seen that $\mathscr{T}$ is JdLG-admissible if and only if $K(\mathscr{S})$ is a (compact topological) group. 
In this case there exists a unique minimal idempotent $Q\in \mathscr{S}$. 
The  \textbf{reversible part} of $E$ is defined by setting $E_{r}:=\operatorname{ran}(Q)$, and as discussed above, it is also described by \eqref{Er-Eaws-definition}. 
Letting $\mathscr{S}_r$ denote the operator semigroup obtained by restriction of the operators in $\mathscr{S}$ to $E_r$, it is known that $\mathscr{S}_r$ is a compact topological group that is (algebraically and topologically) isomorphic to $K(\mathscr{S})$ (cf. \cite[Theorem 16.24]{OTAoET}).
We call a tuple $(e_j)_{j=1}^n\in E^n$ a \textbf{finite unitary system} for $\mathscr{T}$ if the vectors $e_1,...,e_n$ are linearly independent, the set $F:=\operatorname{lin}\{e_1,...,e_n\}$ is $\mathscr{T}$-invariant, and the corresponding matrix representation $\chi$ of $\mathscr{T}$ defined by
    \[
T e_j:=\sum_{i=1}^n \chi_{ij}(T)e_i
\quad (j\in\{1,\ldots,n\},\ T\in\mathscr{T}),
\]
is unitary, i.e. for every $T\in \mathscr{T}$, $\chi(T)$ is a unitary $n\times n$ matrix.
If $(e_j)_{j = 1}^n$ is a finite unitary system for $\mathscr{T}$ with the corresponding matrix representation $\chi$, then for each $1 \le i,j \le n$ the map $T\mapsto \chi_{ij}(T)$ is continuous.
A unitary system $(e_j)_{j=1}^n$ is called \textbf{irreducible} if $F=\operatorname{lin}\{e_1,...,e_n\}$ does not contain any nontrivial $\mathscr{T}$-invariant subspaces. By \cite[Theorem 16.31]{OTAoET} the reversible part of a JdLG-admissible semigroup $\mathscr{T}$ is given by
\begin{equation}\label{CharacterizingReversiblePartEquation}
    E_{r}=\overline{\operatorname{lin}}\{u: (e_j)_{j} \text{ is an irreducible finite unitary system for }\mathscr{T} \text{ and }u=e_j\text{ for some }j\}.
\end{equation}
 
Under an additional assumption one can characterize the reversible part as follows.
\begin{theorem}[cf. \cite{OTAoET}, Proposition 16.29.]\label{relativcompact}
Let $\mathscr{T}$ be a JdLG-admissible semigroup of bounded linear operators on a Banach space $E$. 
If 
\[\inf_{T\in \mathscr{T}}\|Ty\|>0\]
for every $y\neq 0$ such that $\mathscr{T}y$ is relatively compact in $E$, then 
\[E_{r}=\{\xi\in E:\mathscr{T}\xi \text{ is relatively compact in }E\}.\]
    
\end{theorem}
Clearly, Theorem \ref{relativcompact} applies to any JdLG-admissible semigroup of isometries.
Another case of interest is when we have representations of groups.

\begin{corollary}\label{ReversiblePartIsRelativelyCompactOrbits}
    Let $G$ be a group, $E$ be a Banach space and $\pi$ a JdLG-admissible representation of $G$ on $E$. 
    Then 
    \[E_{r}=\{\xi\in E: \pi(G)\xi \text{ is relatively compact in }E\}.\]
    
\end{corollary}
\begin{proof}
Since $\pi$ is JdLG-admissible, it is relatively weakly compact, thus uniformly bounded, so there exists $C>0$ such that $\displaystyle\sup_{g\in G}\|\pi_g\|\leq C$. 
Hence, for any $x\in E$ we obtain 
\[\|x\|=\left\|\pi_g^{-1}\pi_g x\right\|\leq \|\pi_{g^{-1}}\|\|\pi_g x\|\leq C\|\pi_gx\|,\]
and by Theorem \ref{relativcompact} the result follows.
\end{proof}

\begin{example}\label{ex-non-isometric}
    To see that Corollary \ref{ReversiblePartIsRelativelyCompactOrbits} does not hold for general semigroups, let $\pi$ be the representation of $\mathbb{N}$ on $E = \mathbb{C}$ given by $\pi(n)(z) = 2^{-n}z$, and observe that $E = E_{\text{aws}} = \{\xi \in E | \pi(\mathbb{N})\xi\text{ is relatively compact}\}$ and $E_r = \{0\}$.
    Nonetheless, the characterization of $E_r$ given in Equation \eqref{CharacterizingReversiblePartEquation} shows us that every $\xi \in E_r$ has a relatively compact orbit in $E$.
\end{example}

The following proposition presents another particular case of JdLG-admissible semigroups.

\begin{proposition}\label{NewAdditionToJdLGAdmissibleExamples}
Let $E$ be a Banach space and let $\mathscr{S}$ be a weakly compact semigroup of bounded linear operators on $E$ such that every $\mathscr{S}$-orbit contains a unique fixed point. Then $\mathscr{S}$ is JdLG-admissible. Moreover, $E_r$ is the set of fixed vectors of $\mathscr{S}$.
\end{proposition}

\begin{proof}
Let $E_1$ denote the set of fixed points of $\mathscr{S}$, which is a closed linear subspace of $E$.
For any $\xi\in E$ let $\xi_1 = u_1 \xi$ be the fixed point in $\mathscr{S}\xi$. Then $u_1 (\xi-\xi_1) = 0$, so $\xi_0:= \xi-\xi_1\in E_{\text{aws}}$, which shows that $E = E_1 + E_{\text{aws}}$.
The set $E_{\text{aws}}$ is also a linear subspace: if $\xi_1,\xi_2\in E_{\text{aws}}$ and $u_1\xi_1=0$, then $u_1(\xi_1+\xi_2) = u_1\xi_2$. By assumption, $\mathscr{S}u_1\xi_2$ contains a fixed point $vu_1\xi_2$. This is also a fixed point in $\mathscr{S}\xi_2$, so uniqueness implies that $vu_1\xi_2=0$. Hence $vu_1(\xi_1+\xi_2)=0$ and $\xi_1+\xi_2\in E_{\text{aws}}$. (Stability under multiplication by scalars is obvious.) 
The intersection of these subspaces is clearly $0$, so we have a direct sum decomposition.

We next show that $E_{\mathrm{aws}}$ is $\mathscr{S}$-invariant and closed.
For $\xi\in E$, let $P\xi$ denote the unique fixed point in $\mathscr{S}\xi$.
By the direct sum decomposition above, $P$ is linear, 
$\operatorname{ran}(P)=E_1$, and $\ker(P)=E_{\mathrm{aws}}$.
If $\xi\in E_{\mathrm{aws}}$ and $T\in\mathscr{S}$, then the unique fixed
point in $\mathscr{S}(T\xi)$ is also a fixed point in $\mathscr{S}\xi$.
Since the unique fixed point in $\mathscr{S}\xi$ is $0$, we obtain
$T\xi\in E_{\mathrm{aws}}$. Hence $E_{\mathrm{aws}}$ is
$\mathscr{S}$-invariant.

To prove closedness, we show that $P$ has closed graph. Suppose that
$\xi_n\to\xi$ and $P\xi_n\to\eta$. Choose $u_n\in\mathscr{S}$ such that
$u_n\xi_n=P\xi_n$. By weak compactness of $\mathscr{S}$, after passing
to a subnet we may assume that $u_n$ converges in the weak operator
topology to some $u\in\mathscr{S}$. Since $\mathscr{S}$ is uniformly
bounded, $u_n\xi_n\to u\xi$ weakly, and hence $\eta=u\xi$. Since $E_1$
is closed, $\eta$ is fixed, and uniqueness yields $\eta=P\xi$.
Thus $P$ has closed graph and therefore $P\in\mathcal{L}(E)$.
Consequently, $E_{\mathrm{aws}}=\ker(P)$ is closed.

The space $E_1$ is in fact $E_r$ as defined in Theorem \ref{GeneralJdLGDecomposition}: if $\xi\in E_r$ and $\xi=\xi_1+\xi_0$ is its de\-com\-po\-si\-tion as above, then choose $u\in \mathscr{S}$ such that $u\xi_0=0$. For any $v\in \mathscr{S}$ we have then $vu\xi = \xi_1$, so choosing $v$ such that $vu\xi=\xi$ we obtain $\xi=\xi_1$ and thus $E_r\subset E_1$. The converse inclusion is obvious: \eqref{Er-Eaws-definition} holds with any $u,v\in \mathscr{S}$.
\end{proof}

\begin{example}
If the fixed point is not unique, the statement may not hold. Consider the semigroup of two elements $S = \{a,b\}$ with multiplication $a^2 = ba = a$, $b^2 = ab = b$. It acts on $E=C(S)$ by right translations, so that $R_a f \equiv f(a)$ and $R_bf \equiv f(b)$ for every $f$. A function is fixed under this action if and only if it is constant, and is in the set $E_{\text{aws}}$ if it vanishes at $a$ or at $b$. But this set is not a linear subspace, and the decompo\-si\-tion in $E_r+E_{\text{aws}}$, though possible, is not unique.
\end{example}

Proposition \ref{NewAdditionToJdLGAdmissibleExamples} applies to the following special semigroups not covered, in general, by Theorem \ref{GeneralJdLGDecomposition}: 
$\mathscr{S}$ is the weak closure of $\pi(P)$ where $P$ is the semigroup of normal states of a quantum semigroup $M$, and $\pi$ its so-called $\ell_2$-bounded representation on a reflexive Banach space (see the proof of Proposition 3.1 in \cite{kuz-means}). 
One example of such a semigroup $P$ is the set of probability measures on a locally compact group $G$, with convolution as the semigroup operation. 
This particular semigroup is known to have a two-sided invariant mean on $W(P)$ and thus falls under the assumptions of Theorem \ref{GeneralJdLGDecomposition}, but in general this is not the case.


As we mentioned before, every coefficient of a JdLG-admissible representation of a semigroup $S$ is in the space $W(S)$ of weakly almost periodic functions. If $G$ is a group, $W(G)$ contains in particular $B(G)$, the Fourier-Stieltjes algebra of $G$, and as a consequence its uniform closure $\overline{B}(G)$. It is known (see, e.g., \cite{Rudin1959} or \cite{DunklRamirez1971}) that $\overline{B}(G) \subsetneq W(G)$ for any non-compact abelian locally compact group $G$. Moreover, every function in $W(G)$ can be realized as a matrix coefficient of an isometric representation of $G$ on a reflexive Banach space, and such a representation is JdLG-admissible. Therefore, the strict inclusion $\overline{B}(G)\subsetneq W(G)$ also implies that not every JdLG-admissible representation of $G$ is weakly equivalent to a unitary representation on a Hilbert space. In contrast, the following theorem shows that the reversible part of every JdLG-admissible representation is weakly equivalent to a unitary representation.

\begin{theorem}\label{(Weak)EquivalenceForReversiblePart}
Let $S$ be a semigroup, let $E$ be a Banach space, let $\pi$ be a
JdLG-admissible representation of $S$ on $E$, and let $\pi_r$ be
the restriction of $\pi$ to $E_r$. Then $\pi_r$ is weakly equivalent
to a unitary representation $U$ of $S$ on a Hilbert space $\mathcal H$.
If $E$ is a Hilbert space, then $\pi_r$ is equivalent to such a
unitary representation.
\end{theorem}

\begin{proof}
    Let us first consider the case in which $E$ is a Banach space.
    Let $X_n$ denote the set of finite-dimensional  representations $\chi$ of $\pi(S)$ corresponding to some unitary system $(e_j)_{j = 1}^n$ in $E_r$, and let $X = \bigcup_{n = 1}^\infty X_n$.
    For each $n \in \mathbb{N}$ and $\chi \in X_n$, let $\mathcal{H}_\chi = \mathbb{C}^n$, and observe that $\tilde{\chi} := \chi\circ\pi$ can be viewed as a continuous representation of $S$ on $\mathcal{H}_\chi$.
    Let $\mathcal{H} = \displaystyle\oplus_{\chi \in X}\mathcal{H}_\chi$ be the Hilbert space direct sum of these spaces, and observe that $U = \oplus_{\chi \in X}\tilde{\chi}$ is a unitary representation of $S$ on $\mathcal{H}$ that decomposes into a direct sum of finite dimensional representations.
    
    We now show that $U$ is weakly contained in $\pi_r$.
    Every coefficient $\phi(s) = f'(U_s \xi)$ of $U$ is approximated (uniformly on $S$) by coefficients with $\xi$ supported in only finitely many of the spaces $\mathcal H_\chi$, and by linearity, it is sufficient to consider only $\phi$ such that $\xi \in \mathcal{H}_\chi$ with a single $\chi\in X$.
    Fix $\chi \in X$, and let $(e_j)_{j = 1}^n$ be a corresponding unitary system for $\pi_r(S)$. 
    For $F = \text{lin}\{e_1,\cdots,e_n\}$ and for $1 \le i \le n$ we define $f_i^\prime \in F^\prime$ by 
    
    \begin{equation}
        f_i^\prime\left(\sum_{j = 1}^nc_je_j\right) = c_i,
    \end{equation}
    and extend $f_i^\prime$ to all of $E_r$ by the Hahn-Banach theorem. Then we see that $\tilde{\chi}_{ij}(s) = f_i^\prime(\pi(s)e_j)$, which yields the desired result.

    We will now show that $\pi_r$ is weakly contained in $U$.
    Let $(e_j)_{j = 1}^n$ be a finite unitary system in $E_r$ and let $f^\prime \in E^\prime$.
    We see that

    \begin{equation}
        f^\prime\left(\pi(s)e_j\right) = \sum_{i = 1}^n\chi_{ij}(\pi(s))f^\prime(e_i),
    \end{equation}
    so the matrix coefficient $s\mapsto f^\prime(\pi(s)e_j)$ is a linear combination of matrix coefficients of $U$. 
    Equality \eqref{CharacterizingReversiblePartEquation} now yields the desired result.

    Now let us consider the case in which $E$ is a Hilbert space.
    Let $\mathscr{S}$ be the closure of $\pi(S)$ in the weak topology of $\mathcal{L}(E)$, let $\mathscr{S}_r$ denote the restriction of $\mathscr{S}$ to $E_r$, and observe that $\mathscr{S}_r$ is the weak closure of $\pi_r(S)$. 
    We have already seen that $\mathscr{S}_r$ is a compact (hence amenable) topological group.
    Day \cite{day1950means} showed that any uniformly bounded representation of an amenable group $G$ on a Hilbert space is \textit{unitarizable}, i.e., is equivalent to a unitary representation of $G$.
    Applying this to $G=\mathscr{S}_r$, we get that the identity map is equivalent to a unitary representation $U'$ of $\mathscr{S}_r$ on a Hilbert space $\mathcal{H}$.
    The Peter-Weyl Theorem tells us that any unitary representation of a compact group decomposes into a direct sum of finite dimensional irreducible representations.
    The desired result follows from the observation that $\pi_r$ is equivalent to the unitary representation $U=U'\circ \pi_r$.
\end{proof}

\begin{remark}
    Fix some irrational $\alpha \in \mathbb{S}^1$ and consider the representation $\pi$ of $\mathbb{Z}$ on $C(\mathbb{S}^1)$ given by $(\pi_nf)(x) = f(\alpha^nx)$.
    We see that the only finite-dimensional irreducible subrepresentations of $\pi$ are the $1$-dimensional spaces spanned by the eigenvectors $f_n(z) = z^n$.
    However, we cannot view $C(\mathbb{S}^1)$ as the direct sum of subspaces spanned by these monomials.
    This example shows why we need to work with weak equivalence in Theorem \ref{(Weak)EquivalenceForReversiblePart} when dealing with Banach spaces.
\end{remark}

\section{The almost weakly stable part of the JdLG-decomposition}\label{AlmostWeaklyStableSection}
We begin with a characterization of the almost weakly stable part of the JdLG-decomposition through the unique invariant mean on the space of weakly almost periodic functions.
This characterization is related to the flight functions discussed in \cite[Appendix 10]{TempelmanBook}.

\begin{lemma}\label{MeanCharacerizationOfAWS}
    Let $S$ be a semitopological semigroup with a bi-invariant mean $m$ on $W(S)$, let $E$ be a Banach space, and let $\pi$ be a relatively weakly compact representation of $S$ on $E$.
    For $\xi \in E$, the following are equivalent:
    \begin{enumerate}[(i)]
        \item $\xi \in E_{\text{aws}}$.

        \item For all $p > 0$, and all $f^\prime \in E^\prime$, we have
        \begin{equation}\label{m-of-power-p=0}
            m\left(s\mapsto|f'(\pi_s\xi)|^p\right) = 0.
        \end{equation}

        \item There exists $p > 0$ such that for all $f^\prime \in E^\prime$ Equation \eqref{m-of-power-p=0} holds.
    \end{enumerate}
\end{lemma}

\begin{proof}
    We will first show that (i) implies (ii).
    Let $(s_i)_{i \in I} \subseteq S$ be a net for which $\displaystyle\lim_{i}\pi_{s_i}\xi = 0$ weakly, and let $p > 0$ be arbitrary.
    Since, for each $f' \in E'$ the function $s\mapsto |f'(\pi_s\xi)|^p$ belongs to $\text{W}(S)$, we have that $\{|f'(\pi_{ss_i}\xi)|^p\}_{i \in I}$ is a relatively weakly compact set. Hence there exists a subnet $(s_j)_{j \in J}$ for which $|f'(\pi_{ss_j}\xi)|^p$ converges weakly to some function $\phi(s)$.
    Since weak convergence implies pointwise convergence and $\displaystyle\lim_j\pi_{s_j}\xi = 0$ weakly, we see that $\phi$ is the $0$ function, hence

    \begin{equation}
        m(|f'(\pi_s\xi)|^p) = \lim_jm(|f'(\pi_{ss_j}\xi)|^p) = m(\lim_j|f'(\pi_{ss_j}\xi)|^p) = m(0) = 0.
    \end{equation}

    It is clear that (ii) implies (iii), so we will now show that (iii) implies (i).
    
    Let $S^w$ denote the weakly almost periodic compactification of $S$, and we let $i, M,$ and $\mu$ be as in Section \ref{WAPSubsection}. Let
\[
\bar{\pi}:S^w\longrightarrow \overline{\pi(S)}^{\,\mathrm{WOT}}
\]
denote the continuous extension of $\pi$.
Let $f^\prime\in E^\prime$ be arbitrary.
   Since
\begin{equation}
    m\left(\left|f^\prime(\pi_s\xi)\right|^p\right)
    = M\left(i\left(\left|f^\prime(\pi_s\xi)\right|^p\right)\right)
    = \int_{K(S^w)}\left|f^\prime(\bar{\pi}_s\xi)\right|^p\,d\mu(s),
\end{equation} 
we see that $m\left(|f^\prime(\pi_s\xi)|^p\right)=0$ if and only if
$f^\prime(\bar{\pi}_s\xi)=0$ for $\mu$-a.e. $s\in K(S^w)$, and in fact
for every $s\in K(S^w)$ since this is a continuous function.
Since $f^\prime$ is arbitrary, this implies
$\bar{\pi}_s\xi=0$ for every $s\in K(S^w)$.
By the definition in \eqref{Er-Eaws-definition}, we have
$\xi\in E_{\mathrm{aws}}$.
\end{proof}

Our next result can be seen as a generalization of \cite{jones1976ergodic}, \cite[Theorem 16.34]{OTAoET}, and \cite[Theorem 2.5]{eisner2007weakly}.
Furthermore, it relates the almost weakly stable vectors to the vectors from $\mathcal{H}_w$ in the compact-weak mixing decomposition. 

\begin{theorem}\label{EquivalentCharacterizationsOfAWSWhenAmenable}
    Let $(S,\lambda)$ be a bi-amenable measured semigroup.
    Let $E$ be a Banach space, and let $\pi$ be a relatively weakly compact (hence JdLG-admissible) representation of $S$ on $E$.
    Given $\xi \in E$ and a left-invariant mean $m$ on $L^\infty(S,\lambda)$, (i)-(iv) are equivalent.

    \begin{enumerate}[(i)]
        \item $\xi \in E_{\text{aws}}$, i.e., $0 \in c\ell_\sigma(\pi(S)\xi)$.

        \item For some/each $p > 0$, $\displaystyle m\left(s\mapsto\left|\left\langle \pi_s\xi, x^\prime\right\rangle\right|^p\right) = 0\ \forall\ x^\prime \in E'$.
        
        \item For any weakly open neighborhood $U$ of $0 \in E$, $m\left(\left\{s \in S | \pi_s\xi \in U\right\}\right) = 1$.
        
        \item $d_*-\displaystyle\lim_s \pi_s\xi = 0\text{ weakly}$.
    \end{enumerate}
    \noindent Moreover, if $S$ is locally compact, $\lambda$ is a Radon measure, and $(S,\lambda)$ possesses a left-F\o lner net $\mathcal{F} = (F_i)_{i \in I}$, then (i)-(viii) are equivalent.
    \begin{enumerate}
        \item[(v)] For some/each $p > 0$, $\displaystyle\lim_i\frac{1}{\lambda(F_i)}\int_{F_i}\left|\left\langle \pi_s\xi, x^\prime\right\rangle\right|^pd\lambda(s) = 0\ \forall\ x^\prime \in E'$.
        
        \item[(vi)] $d_{\mathcal{F}}-\displaystyle\lim_s \pi_s\xi = 0\text{ weakly}$.
        
        \item[(vii)]$\underline{d}-\displaystyle\lim_s \pi_s\xi = 0\text{ weakly}$.
    \end{enumerate}
    \begin{enumerate}
        \item[(viii)] For some/each $p > 0$, $\displaystyle\lim_{i}\sup_{\|x^\prime\|\leq 1}\frac{1}{\lambda(F_i)}\int_{F_i}\left|\left\langle \pi_s\xi, x^\prime\right\rangle\right|^p d\lambda(s) = 0$.
    \end{enumerate}
    \noindent \noindent If we further assume $\mathcal{F} = (F_n)_{n = 1}^\infty$ is a F\o lner sequence and that $\displaystyle\lim_{n\rightarrow\infty}\lambda(F_n) = \infty$, then (i)-(ix) are equivalent.
    \begin{enumerate}
        \item[(ix)] There exists a set $A \subseteq S$ with $d_\mathcal{F}(A) = 1$ such that $\displaystyle \lim_{\underset{s \in A}{s\to \infty}}\pi_s\xi = 0$ weakly. 
    \end{enumerate}
\end{theorem}

\begin{proof}
    As discussed in Section \ref{WAPSubsection}, $m$ restricts to a bi-invariant mean on $W(S)$.
    We now use Lemma \ref{MeanCharacerizationOfAWS} to see that for some/each $p > 0$ we have
    \begin{alignat*}{2}
        E_{\text{aws}} 
        &=\left\{\xi \in E\ :\ m\left(s\mapsto\left|\left\langle \pi_s\xi, x^\prime\right\rangle\right|^p\right) = 0\ \forall\ x' \in E'\right\},     
    \end{alignat*}
    hence (i) and (ii) are equivalent.

    To see that (ii) implies (iii), let $U$ be a weak neighborhood of $0$.
Then there exist $x_1^\prime,\ldots,x_n^\prime\in E^\prime$ and
$\epsilon>0$ such that
\[
V:=\bigcap_{j=1}^n
\left\{\eta\in E:\left|\langle\eta,x_j^\prime\rangle\right|^p<\epsilon\right\}
\subseteq U.
\]
For each $j=1,\ldots,n$, we have
\[
0=m\left(s\mapsto
\left|\left\langle\pi_s\xi,x_j^\prime\right\rangle\right|^p\right)
\geq
\epsilon\,m\left(\mathbbm{1}_{V_j^c}(\pi_s\xi)\right),
\]
where
\[
V_j:=\left\{\eta\in E:
\left|\langle\eta,x_j^\prime\rangle\right|^p<\epsilon\right\}.
\]
Thus $m(\mathbbm{1}_{V_j^c}(\pi_s\xi))=0$ for every $j$. Since
$V^c=\bigcup_{j=1}^n V_j^c$, it follows that
$m(\mathbbm{1}_{V^c}(\pi_s\xi))=0$, and hence
\[
m(\mathbbm{1}_U(\pi_s\xi))=1.
\]
   To see that (iii) implies (ii), let $x^\prime\in E^\prime$ and
$\epsilon>0$, and put $C:=\sup_{s\in S}
\left|\left\langle\pi_s\xi,x^\prime\right\rangle\right|^p<\infty.$
Let
\[
U:=\left\{\eta\in E:
\left|\langle\eta,x^\prime\rangle\right|^p<\epsilon\right\}.
\]
By (iii), we have
$m(\mathbbm{1}_{U^c}(\pi_s\xi))=0$. Hence
\[
m\left(s\mapsto
\left|\left\langle\pi_s\xi,x^\prime\right\rangle\right|^p\right)
\leq
\epsilon\,m\left(\mathbbm{1}_U(\pi_s\xi)\right)
+C\,m\left(\mathbbm{1}_{U^c}(\pi_s\xi)\right)
\leq\epsilon.
\]
Since $\epsilon>0$ was arbitrary, (ii) follows.

    It is clear that (iv) implies (iii).
To see that (iii) implies (iv), recall that (iii) implies (ii), hence (i).
Since (i) is independent of the choice of the left-invariant mean, it follows that
(iii) holds for every left-invariant mean. Therefore, every left-invariant mean
assigns value $1$ to the relevant return set, and hence its lower Banach density
is $1$.

    Now let us assume that $(S,\lambda)$ admits a left-F\o lner net $\mathcal{F} = (F_i)_{i \in I}$.
    We observe that $\phi(s) := |\langle \pi_s\xi,x^\prime\rangle|^p$ is in $W(S)$.
    For each $m \in \mathcal{M}(\mathcal{F})$, $m$ restricts to the unique bi-invariant mean on $W(S)$, so the value of $m(\phi)$ is independent of the choice of $m$, hence
    \begin{equation}\label{AveragesExistEquation}
        \lim_i\frac{1}{\lambda(F_i)}\int_{F_i}\left|\left\langle\pi_s\xi,x^\prime\right\rangle\right|^pd\lambda(s) = m(\phi)\text{ for all }m \in \mathcal{M}(\mathcal{F}).
    \end{equation}
    In particular, the limit on the left-hand side of \eqref{AveragesExistEquation} exists.
    We recall that (ii) is equivalent to (i), and (i) is independent of the choice of $m$.
    In particular, if (ii) holds for a particular left-invariant mean, then (ii) holds for all left-invariant means, which shows that (ii) is equivalent to (v).
    Similarly, we see that (iii) is equivalent to (vi).

    It is clear that (vii) implies (vi).
   To see that (vi) implies (vii), recall that (vi) is equivalent to (i).
Since (i) is independent of the choice of the left-F\o lner net, 
applying (i) implies (vi) to every left-F\o lner net $\mathcal{F}$ 
shows that
\[
d_{\mathcal{F}}(\{s\in S:\pi_s\xi\in U\})=1
\]
for every weak neighborhood $U$ of zero and every left-F\o lner net 
$\mathcal{F}$. Taking the infimum over all left-F\o lner nets yields
\[
\underline{d}(\{s\in S:\pi_s\xi\in U\})=1,
\]
which proves (vii). 
   It is now clear that (viii) implies (v), so we will proceed to show that (v) implies (viii).
    Since $(\pi_s)_{s\in S}$ is uniformly bounded, we suppose without loss of generality (as mentioned in Section \ref{RepresentationsOfSemigroupsSubsection}) that $S$ is a monoid and that $(\pi_s)_{s\in S}$ consists of contractions.
    Then the adjoint action $(\pi_s^*)_{s\in S}$ leaves the ($\sigma^*$-compact) dual unit ball $D^\prime \subseteq E^\prime$ invariant. 
    Let $(T_s)_{s\in S}$ denote the Koopman action on $C(D^\prime)$ induced by $(\pi^*_s)_{s\in S}$. 
    For $\xi$ satisfying (ii) we let $f_\xi\in C(D^\prime)$, $f_\xi(x^\prime):=\left|\langle  \xi,x^\prime \rangle\right|^p$, $p>0$. 
    
    We will first show that the convex hull of $T(S)f_\xi$ is relatively weakly compact in $C(D^\prime)$.
    To this end, let $(s_k)_{k\in\N}$ be any sequence in $S$.
    Since $\pi$ is a relatively weakly compact representation, $\pi(S)\xi$ is relatively weakly compact, so the Eberlein-\v{S}mulian Theorem tells us that $\text{cl}_\sigma(\pi(S)\xi)$ is weakly sequentially compact.
    Consequently, there exists a subsequence $(s_{k_j})_{j\in\N}$ and $\eta \in E$ such that \[ \lim_j T_{s_{k_j}}f_\xi(x^\prime)=\lim_j |\langle\pi_{s_{k_j}}\xi,x^\prime \rangle|^p=|\langle\eta,x^\prime \rangle|^p=f_\eta(x^\prime)\] holds 
    for all $x^\prime \in D^\prime$. 
    By the dominated convergence theorem we obtain $\lim_j T_{s_{k_j}}f_\xi=f_\eta$ weakly in $C(D^\prime)$.
    We invoke the Eberlein-\v{S}mulian Theorem once again to see that $T(S)f_\xi$ is relatively weakly compact, from which the desired result follows.
    
    For a F\o lner net $(F_i)_{i\in I}$ we set  \[A_if_\xi(x^\prime):=\displaystyle\frac{1}{\lambda(F_i)}\int_{F_i}T_sf_\xi(x^\prime) d\lambda(s) \]
    for any $x^\prime\in D^\prime$ and $i\in I$. Since the net $(A_i f_\xi)_{i\in I}$ lies in the weak closure of the convex hull of $\{T_s f_\xi : s\in S\}$, it admits a cluster point. 
    By assumption (v) we have $\lim_i A_if_{\xi}(x^\prime)=0$ for any $x^\prime\in D^\prime$, hence every cluster point of $(A_if_\xi)_{i\in I}$ must be the zero function. 
    Eberlein's generalized mean ergodic theorem (see, e.g., Chapter 2, Theorem 1.5 of \cite{UKErgodicThms}) implies
\[0=\lim_{i}\left\|\frac{1}{\lambda(F_i)}\int_{F_i}T_sf_\xi d\lambda(s)\right\|_\infty=\lim_{i}\sup_{\|x^\prime\|\leq 1}\frac{1}{\lambda(F_i)}\int_{F_i}|\langle \pi_s\xi, x^\prime\rangle|^p d\lambda(s). \]

    Now let us show that (vi) is equivalent to (ix) under the additionally stated assumptions.
    Since $d_\mathcal{F}(C) = 0$ for every compact set $C$, we see that (ix) implies (vi), so it remains to show that (vi) implies (ix).
    As in the previous case, we can suppose that every $\pi_s$ is a contraction. 
    Define $\phi(s) = \pi_s\xi$ and observe that $\phi:S\rightarrow E$ is weakly continuous, so $K := \text{cl}_\sigma(\phi(S))$ is weakly compact.
    For each $n \in \mathbb{N}$, let $\nu_n$ be the probability measure on $K$ given by 

\begin{equation}
    \nu_n(B)
    = \frac{\lambda\left(F_n\cap\phi^{-1}(B)\right)}
    {\lambda(F_n)}.
\end{equation}    
Since $\lambda$ is a Radon measure, $\lambda|_{F_n}/\lambda(F_n)$ is a Radon probability measure, and $\nu_n$ is its pushforward under the continuous map $\phi$. Hence $\nu_n$ is a Radon probability measure on $K$.    It follows that $\nu := \sum_{n = 1}^\infty2^{-n}\nu_n$ is also a Radon probability measure on $K$, so $\nu(K\setminus\text{supp}(\nu)) = 0$.
    Let $H = \text{supp}(\nu)$, and observe that $\nu(U) > 0$ for every relatively open set $U$ in $H$, so $H$ satisfies the countable chain condition.
    In other words, any collection of disjoint open sets in $H$ must be countable.
    The main result of \cite{CCCEberleinCompactaAreMetrizable} tells us that a weakly compact subset of a Banach space is metrizable if and only if it satisfies the countable chain condition, hence $H$ is metrizable.
    Since (vi) holds, we see that $0 \in H$, so let $\{U_n\}_{n = 1}^\infty$ be a decreasing countable neighborhood basis for $0$.
    Set $B:=\phi^{-1}(H)$. Since $\nu_n(K\setminus H)=0$, we have
\[
\lambda(F_n\setminus B)=0
\]
for every $n\in\mathbb{N}$.
Choose the decreasing neighborhood basis $(U_n)_{n\in\mathbb{N}}$ of $0$ in $H$ and write
$U_n=V_n\cap H$, where $V_n$ is weakly open in $E$.
For each $n\in\mathbb{N}$, set
\[
A_n:=B\cap\phi^{-1}(V_n)=B\cap\phi^{-1}(U_n).
\]
Then $(A_n)_{n\in\mathbb{N}}$ is decreasing and
$d_{\mathcal{F}}(A_n)=1$ for every $n\in\mathbb{N}$.
    Let $m_n$ be such that $\lambda(A_n\cap F_m) > \frac{n-1}{n}\lambda(F_m)$ for all $m \ge m_n$.
    Assuming without loss of generality that $m_{n+1} > m_n$, we see that for $A := \bigcup_{n = 1}^\infty\bigcup_{m = m_n+1}^{m_{n+1}}(A_n\cap F_m)$ we have $d_\mathcal{F}(A)=1$ and
    $$
    \{ s\in A: \pi_s \xi\notin U_n\} = A\setminus A_n \subset \bigcup_{k<n}\;\bigcup_{m = m_k+1}^{m_{k+1}}(A_k\cap F_m) \subseteq \bigcup_{m = 1}^{m_n}F_m.
    $$
\end{proof}

\begin{corollary}\label{cor:ConvergenceAlongDensity1Set}
    Let $(S,\lambda)$ be a locally compact bi-amenable measured semigroup with $\lambda$ a Radon measure, $E$ a Banach space, and $\pi$ a relatively weakly compact representation of $S$ on $E$.
    Suppose that $(S,\lambda)$ admits a left-F\o lner sequence $\mathcal{F} = (F_n)_{n = 1}^\infty$ with $
\lim_{n\to\infty}\lambda(F_n)=\infty$, and that $E_{\text{aws}}$ is norm separable.
    There exists an $A \subseteq S$ with $d_\mathcal{F}(A) = 1$ such that

    \begin{equation}\label{eq:SingleSequenceForAllOfE_aws}
        E_{\text{aws}} = \left\{\xi \in E\ |\ \lim_{\underset{s \in A}{s\rightarrow\infty}}\pi_s\xi = 0\text{ weakly}\right\}.
    \end{equation}
\end{corollary}

\begin{proof}
    The equivalence of parts (vi) and (ix) of Theorem \ref{EquivalentCharacterizationsOfAWSWhenAmenable} shows us that the right hand side of Equation \eqref{eq:SingleSequenceForAllOfE_aws} is contained in the left hand side, so it only remains to show the reverse inclusion.
    To this end, let $(\xi_n)_{n = 1}^\infty$ be dense in $E_{\text{aws}}$.
    For each $n$ let $A_n \subseteq S$ be such that $d_\mathcal{F}(A_n) = 1$ and $\displaystyle\lim_{\underset{s \in A_n}{s\rightarrow\infty}}\pi_s\xi_n = 0$ weakly.
    Let $B_n = \bigcap_{m = 1}^nA_m$, and observe that $d_\mathcal{F}(B_n) = 1$ for all $n$.
    Let $m_n$ be such that $\lambda(B_n\cap F_m) \ge \frac{n-1}{n}\lambda(F_m)$ for all $m \ge m_n$ and $m_{n+1} > m_n$.
    The desired set $A$ is given by $A := \bigcup_{n = 1}^\infty\bigcup_{m = m_n+1}^{m_{n+1}}(B_n\cap F_m)$, where $m_0 = 0$.
\end{proof}

\begin{example}\label{CounterExampleForAWSCharacterizations}
    We now recall \cite[Example 4.3]{klawe1980dimensions} to show that Lemma \ref{MeanCharacerizationOfAWS} does not necessarily hold if $W(S)$ only admits a left-invariant mean and not a bi-invariant mean.
    This example will also show that Theorem \ref{EquivalentCharacterizationsOfAWSWhenAmenable} does not necessarily hold if $S$ is a left-amenable subsemigroup of an amenable group that is not right-amenable.

    We consider the semigroup $(S,\cdot)$ given by $S = \{(m,n) \in \mathbb{Z}^2\ :\ m \ge 0,n\ge 1\}$ with the multiplication
    $$
    (m_1,n_1)\cdot(m_2,n_2) = (2^{n_2}m_1+m_2,n_1+n_2).
    $$
    As shown in \cite{klawe1980dimensions}, $S$ is cancellative and left-amenable, hence it satisfies the strong F\o lner condition, and since it is countable, it possesses a left-F\o lner sequence $\mathcal{F} = (F_i)_{i = 1}^\infty$. 
    Furthermore, by replacing $F_n$ with $F_n\cdot(0,1)$, which still form a left F\o lner sequence, we may assume without loss of generality that $F_n \subseteq S\cdot(0,1)$ for all $n$.
    We observe that
    $$
    S\cdot(0,1) = \{(m,n) \in S\ :\ m\in 2\Z,\; n\ge 2\},
    \quad S\cdot(1,1) = \{(m,n) \in S\ :\ m\notin 2\Z,\; n\ge 2\}.
    $$
    Since $S\cdot(0,1)\cap S\cdot(1,1) = \emptyset$, we see that $S$ is not right-amenable.
    Now let $\xi:S\rightarrow\{0,1\}$ be the characteristic function of $S\cdot (0,1)$. We note that for $t\in S$ and $s=(m,n)\in S$, we have $ts\in S\cdot(0,1)$ if and only if $m\in 2\Z$.
    We see that $\{R_s\xi\}_{s \in S} = \{0,1\}$,
    hence it is relatively weakly compact and $\xi \in W(S)$.
    Now set $E = W(S)$, $\pi = R$, and let $x^\prime$ denote the evaluation at $s_0$ for some $s_0 \in S$.
    We see that $\xi \in E_{\text{aws}}$ because $R_{(1,1)}\xi = 0$.
    However, for every $s\in F_i\subset S\cdot(0,1)$ we have $R_s\xi=1$, so if $m \in \mathcal{M}(\mathcal{F})$, then 

    \begin{equation}
        \lim_i\frac{1}{\lambda(F_i)}\int_{F_i}\left|\left\langle \pi_s\xi,x^\prime\right\rangle\right|d\lambda(s) = m(s\mapsto \left(R_s\xi\right)(s_0)) = m(1) = 1.
    \end{equation}
\end{example} 

\begin{example}
    We now give an example showing that parts (vi) and (ix) of Theorem \ref{EquivalentCharacterizationsOfAWSWhenAmenable} need not be equivalent if $\mathcal{F}$ is not a F\o lner sequence.

    Consider the discrete commutative semigroup $S = \mathbb{Z}\times\mathbb{N}_0$ with addition given by $(a,m)+(b,n) = (a+b,\text{max}(m,n))$.
    We order $\mathbb{N}\times\mathbb{N}_0$ by $(a,b) \ge (c,d)$ if $a \ge c$ and $b \ge d$.
    Let $F_{L,n} = \{1,\cdots,L\}\times\{n\}$ and observe that $\mathcal{F} := \{F_{L,n}\ |\ (L,n) \in \mathbb{N}\times\mathbb{N}_0\}$ is a F\o lner net in $S$.
    Now consider the unitary operator $U$ on the Hilbert space $E = \ell^2(\mathbb{Z})$ given by $Ue_k = e_{k+1}$, where $(e_k)_{k \in \mathbb{Z}}$ is the standard basis for $E$.
    Define a representation $\pi$ of $S$ on $E$ by $\pi_{(a,n)} = U^a$, let $\xi = e_0$, and observe that $\pi_{(a,n)}\xi = e_a$.
    The representation $\pi$ is JdLG-admissible because it consists of contractions on a Hilbert space, and we see that $\xi \in E_{\text{aws}}$.
    We still have that parts (i) and (vi) of Theorem \ref{EquivalentCharacterizationsOfAWSWhenAmenable} are equivalent, so it only remains to show that part (ix) does not hold.

    To this end, let $A \subseteq S$ be such that

    \begin{equation}
        \lim_{\underset{s \in A}{s\rightarrow\infty}}\pi_s\xi = 0\text{ weakly}.
    \end{equation}
    For $a \in \mathbb{Z}$, let $A_a := \{n \in \mathbb{N}_0\ |\ (a,n) \in A\}$.
    We see that for any $a \in \mathbb{Z}$ and any sequence $(n_k)_{k = 1}^\infty$ of pairwise distinct elements of $A_a$, we have

    \begin{equation}
        \lim_{k\rightarrow\infty}\langle\pi_{(a,n_k)}\xi,e_a\rangle = \langle e_a,e_a\rangle = 1 \neq 0.
    \end{equation}
    It follows that we cannot have $(a,n_k)\rightarrow\infty$ in the discrete space $S$, so $A_a$ must be a finite set.
    Now let $(L_0,n_0) \in \mathbb{N}\times\mathbb{N}_0$ be arbitrary and pick $n > n_0$ for which $(1,n),(2,n),\dots,(L_0,n) \not\in A$, i.e., we pick $n > n_0$ such that $F_{L_0,n} \cap A = \emptyset$.
    It follows that arbitrarily far out in the F\o lner net $\mathcal{F}$ we have

    \begin{equation}
        \frac{|A\cap F_{L_0,n}|}{|F_{L_0,n}|} = 0,
    \end{equation}
    so we cannot have $d_\mathcal{F}(A) = 1$.
\end{example}

\begin{example}
    We now give an example to show that Theorem \ref{EquivalentCharacterizationsOfAWSWhenAmenable}(ix) cannot be strengthened to further require that $\underline{d}(A) = 1$.

    It is well known that there exists a weakly mixing rigid unitary operator $U$ on a Hilbert space $\mathcal{H}$ (see, e.g., \cite[Theorem $\mathrm{IV}.3.11$]{TanjaEisnerStabilityBook}).
    More concretely, let $U$ be such that $\mathcal{H} = \mathcal{H}_w$ in the compact-weak mixing decomposition, and there exists a sequence $(n_k)_{k = 1}^\infty$ for which $\displaystyle\lim_{k\rightarrow\infty}U^{n_k} = \text{Id}$, with convergence taking place in the strong operator topology.
    Fix $0\neq \xi \in \mathcal{H}$, and let us assume for the sake of contradiction that there exists $A \subseteq \mathbb{Z}$ with $\underline{d}(A) = 1$ and $\displaystyle\lim_{\underset{n \in A}{n\rightarrow\infty}}U^n\xi = 0$.\footnote{We warn the reader that in other texts the density functions $\overline{d}$ and $\underline{d}$ have different definitions on $\mathbb{Z}$ than what we introduced in Section \ref{AmenableSemigroupsSubsection}. Nonetheless, we are using the definitions from Section \ref{AmenableSemigroupsSubsection}}
    Since $\overline{d}(A^c) = 0$, $A^c$ cannot contain arbitrarily long intervals, so let $L \in \mathbb{N}$ be such that $A\cap[x,x+L] \neq \emptyset$ for all $x \in \mathbb{Z}$.
    Let $\eta \in \mathcal{H}$ be such that $\langle U^n\xi,\eta\rangle \neq 0$ for $0 \le n \le L$.
    Let $0 \le \ell \le L$ be such that there is an infinite subsequence $(m_k)_{k = 1}^\infty$ of $(n_k)_{k = 1}^\infty$ for which $m_k+\ell \in A$ for all $k$.
    The desired contradiction now follows from the fact that $\displaystyle\lim_{k\rightarrow\infty}\langle U^{m_k+\ell}\xi,\eta\rangle = \langle U^\ell\xi,\eta\rangle \neq 0$.
\end{example}

\begin{example}
    We now give an example to show that Corollary \ref{cor:ConvergenceAlongDensity1Set} does not necessarily hold if we do not assume that $E_{\text{aws}}$ is norm separable.

    Let $\mathcal{R}$ denote the collection of increasing sequences of natural numbers $(n_m)_{m = 1}^\infty$ for which we also have $\lim_{m\rightarrow\infty}\frac{n_{m+1}}{n_m} = \infty$.
    Using \cite[Proposition 3.15]{RigidityAndNonRecurrence}, we see that for each $r:= (n_m)_{m = 1}^\infty \in \mathcal{R}$, there exists a weakly mixing unitary operator $U_r$ on a Hilbert space $\mathcal{H}_r$ so that $U_r$ is rigid along $r$, i.e., $\lim_{m\rightarrow\infty}U_r^{n_m} = \text{Id}$ strongly.
    We see that if $A \subseteq \mathbb{N}$ is infinite, then there exists some $r := (n_m)_{m = 1}^\infty \in \mathcal{R}$ contained in $A$.
    Let $0 \neq \xi \in \mathcal{H}_r$ be arbitrary and observe that $\lim_{m\rightarrow\infty}U_r^{n_m}\xi = \xi$, so $\lim_{\underset{n\rightarrow\infty}{n \in A}}U_r^n\xi \neq 0$ weakly.
    Consequently, it suffices to take $U = \oplus_{r \in \mathcal{R}}U_r$ and $\mathcal{H} = \oplus_{r \in \mathcal{R}}\mathcal{H}_r$.
\end{example}

\section{The JdLG decomposition via ultrafilters}\label{UltrafilterSection}

Let $E = \mathcal{H}$ be a Hilbert space, $U$ a unitary operator on $\mathcal{H}$, and $\mathcal{H} = \mathcal{H}_c\oplus\mathcal{H}_w$ the compact-weak mixing decomposition.
In \cite[Corollary 4.6]{MinimalIdempotentsAndERT} Bergelson showed that the orthogonal projection $Q:\mathcal{H}\rightarrow\mathcal{H}_c$ is given by $p-\lim_nU^n$, where $p \in \beta\mathbb{Z}$ is any minimal idempotent.\footnote{In fact, Bergelson proved his result for unitary representations of any countably infinite discrete group.}
Variations of this result appear throughout the literature, such as \cite[Theorem 2.2]{bergelson2007central}, \cite[Theorem 2.25]{BMIPUltrafiltersAndSzemerediForGPolynomials}, \cite[Theorem 4]{mccutcheon2014d}, and \cite[Section 6]{bergelson2021iterated}, and they are used to study the largeness of return times in dynamical systems.

Theorem \ref{ProjectionsViaUltrafilters} generalizes Bergelson's result to JdLG-admissible representations of semigroups on a Banach space $E$.
Theorem \ref{thm:CharacterizationsViaVisits} characterizes the reversible part using largeness of return times, and gives a related characterization for the almost weakly stable part.
For the rest of this section, we freely use the notation and facts from Subsection \ref{UltrafiltersSubsection}.

\begin{lemma}\label{UltrafilterLemma}
    Let $S$ be a discrete semigroup, let $E$ be a Banach space, let $\pi$ be a JdLG-admissible representation of $S$ on $E$, and let $E = E_r\oplus E_{\text{aws}}$ be the JdLG decomposition.
    \begin{enumerate}[(i)]
        \item If $p \in \beta S$ is an idempotent, then $\displaystyle p-\lim_s\pi_s\xi_r = \xi_r$ for all $\xi_r \in E_r$.

        \item If $S$ is cancellative and left reversible, $p\in\beta S$, and
$p^{-1}p\in\beta S$, then $p^{-1}p-\lim_s\pi_s\xi_r=\xi_r$ for all $\xi_r\in E_r$.

        \item If $p \in \overline{K(\beta S)}$, then $p-\displaystyle\lim_s\pi_s\xi_{\text{aws}} = 0$ for all $\xi_{\text{aws}} \in E_{\text{aws}}$.

        \item If $S$ is a bi-amenable semigroup and $p \in \mathcal{D}^*(S)$, then $\displaystyle p-\lim_s\pi_s\xi_{\text{aws}} = 0$ for all $\xi_{\text{aws}} \in E_{\text{aws}}$.
    \end{enumerate}
\end{lemma}

\begin{proof}
    As usual, we let $\mathscr{S} = c\ell_\sigma(\pi(S))$ and we recall that $\mathscr{S}$ is a compact semitopological semigroup.
    We let $\mathscr{S}_r$ denote the restriction of $\mathscr{S}$ to $E_r$, and we recall that $\mathscr{S}_r$ is isomorphic to $K(\mathscr{S})$ as a topological group.
    
    We first prove (i).
    Let $h:\beta S\rightarrow\mathscr{S}_r$ denote the map $h(q) = \displaystyle q-\lim_s\pi_s$ with convergence taking place in the weak operator topology of $\mathcal{L}(E_r)$.
    Since $\pi:S\rightarrow\mathscr{S}_r$ is a homomorphism with dense image in a semitopological semigroup and $\pi = h\circ e$ on $S$, we see that $h$ is a continuous surjective homomorphism.
    In particular, $h$ sends idempotents to idempotents.
    Since $\mathscr{S}_r$ is a compact topological group, it contains a unique idempotent $Q$.
    Furthermore, we have already seen that $Q$ acts as the identity on $E_r$, so the desired result follows.

We now proceed to prove (ii).
Let $G$ denote the group of right quotients of $S$, and let
\[
h_0:S\to\mathscr{S}_r,\qquad h_0(s)=\pi_s|_{E_r}.
\]
Since $h_0$ is a homomorphism into the group $\mathscr{S}_r$,
Proposition \ref{prop:ExtendingSemigroupHomomorphisms} yields a
homomorphism $
\widetilde h:G\to\mathscr{S}_r$ extending $h_0$. Since $\mathscr{S}_r$ is a compact topological group,
$\widetilde h$ extends continuously to
$\widehat h:\beta G\to\mathscr{S}_r.$
Thus, for $p\in\beta S\subseteq\beta G$ we have
\[
\widehat h(p^{-1}p)
=\widehat h(p)^{-1}\widehat h(p)
=Q,
\]
where $Q$ is the identity element of $\mathscr{S}_r$.
If $p^{-1}p\in\beta S$, this means precisely that
\[
p^{-1}p-\lim_s\pi_s|_{E_r}=Q.
\]
Since $Q$ acts as the identity on $E_r$, we obtain $p^{-1}p-\lim_s\pi_s\xi_r=\xi_r$ for every $\xi_r\in E_r$.
    
  We now proceed to prove (iii).
Let $\phi:\beta S\rightarrow\mathscr{S}$ be given by
\[
\phi(q)=q-\lim_s\pi_s.
\]
As in part (i), $\phi$ is a continuous surjective homomorphism.
Hence $
K(\mathscr{S})=\phi(K(\beta S)).$ Since $K(\mathscr{S})$ is closed, it follows that
$\phi\left(\overline{K(\beta S)}\right)=K(\mathscr{S}).$ Now let $p\in\overline{K(\beta S)}$ and set $u:=\phi(p)\in K(\mathscr{S})$.
Let $Q$ denote the identity element of the compact group $K(\mathscr{S})$.
Then $u=uQ$. Hence, for every $\xi_{\mathrm{aws}}\in E_{\mathrm{aws}}=\ker Q$ we obtain
\[
p-\lim_s\pi_s\xi_{\mathrm{aws}}
=
u\xi_{\mathrm{aws}}
=
uQ\xi_{\mathrm{aws}}
=
0.
\]

    We now proceed to prove (iv). 
    Theorem \ref{EquivalentCharacterizationsOfAWSWhenAmenable}(iv) tells us that for any weakly open neighborhood $U$ of $0$ and every $\xi_{\text{aws}} \in E_{\text{aws}}$, we have

    \begin{equation}
        d^*\left(\left\{s \in S\ :\ \pi_s\xi_{\text{aws}} \notin U\right\}\right) = 0\text{, hence }\left\{s \in S\ :\ \pi_s\xi_{\text{aws}} \notin U\right\} \notin p.
    \end{equation}
    Since $p$ is an ultrafilter, we see that $\left\{s \in S\ :\ \pi_s\xi_{\text{aws}} \in U\right\} \in p$.
    Since $U$ was arbitrary, we see that $\displaystyle p-\lim_s\pi_s\xi_{\text{aws}} = 0$.
\end{proof}

\begin{theorem}\label{ProjectionsViaUltrafilters}
    Let $S$ be a discrete semigroup, let $E$ be a Banach space, and let $\pi$ be a JdLG-admissible representation of $S$ on $E$.
    \begin{enumerate}[(i)]
        \item If $p \in \overline{K(\beta S)}$ is an idempotent, then $\displaystyle \pi_p := p-\lim_s\pi_s = Q$, where $Q$ is the projection onto $E_r$ whose kernel is $E_{\text{aws}}$.

        \item If $S$ is bi-amenable, then for each essential idempotent $p \in \mathcal{D}^*(S)$ we have $\pi_p = Q$.

        \item If $S$ is cancellative and left reversible, and $p \in \beta S$ is of the form $q^{-1}q$ for some $q \in \overline{K(\beta S)}$, then $\pi_p = Q$.

      \item If $S$ is cancellative and bi-amenable, and $p \in \beta S$
is of the form $q^{-1}q$ for some $q \in \mathcal{D}^*(S)$,
then $\pi_p=Q$.
    \end{enumerate}
\end{theorem}

\begin{proof}
    Let $\xi \in E$ be arbitrary and let $\xi_r \in E_r$ and $\xi_{\text{aws}} \in E_{\text{aws}}$ be such that $\xi = \xi_r+\xi_{\text{aws}}$.
    Part (i) follows from parts (i) and (iii) of Lemma \ref{UltrafilterLemma}.
    Part (ii) follows from parts (i) and (iv) of Lemma \ref{UltrafilterLemma}.
   To prove part (iii), we recall that $\overline{K(\beta S)}$ is a left ideal in $\beta G$.
Hence, for $q\in\overline{K(\beta S)}$, we have
$q^{-1}q\in\overline{K(\beta S)}$.
The desired result now follows from parts (ii) and (iii) of Lemma
\ref{UltrafilterLemma}. To prove part (iv), we recall that $\mathcal{D}^*(S)$ is a left ideal in $\beta G$.
Hence, for $q\in\mathcal{D}^*(S)$, we have
$q^{-1}q\in\mathcal{D}^*(S)$.
The desired result now follows from parts (ii) and (iv) of Lemma
\ref{UltrafilterLemma}.
\end{proof}

We remind the reader that the notions of largeness used in the next result are discussed in Section \ref{sssec:NotionsOfLargeness}.

\begin{theorem}\label{thm:CharacterizationsViaVisits}
    Let $S$ be a discrete semigroup, let $E$ be a Banach space, and let $\pi$ be a JdLG-admissible representation of $S$ on $E$.
    \begin{enumerate}[(i)]
        \item We have $\xi \in E_r$ if and only if for any $\epsilon > 0$, the set $\{s \in S\ |\ \|\xi-\pi_s\xi\| < \epsilon\}$ is IP$^*$.

        \item We have $\xi \in E_{\text{aws}}$ if and only if for any weakly open neighborhood $U$ of $0$, the set $\{s \in S\ |\ \pi_s\xi \in U\}$ is thickly syndetic.

        \item Assume that $S$ is cancellative and left reversible and let $\ell \in \mathbb{N}$ be arbitrary.
        We have $\xi \in E_r$ if and only if for any $\epsilon > 0$ the set $\{s \in S\ |\ \|\xi-\pi_s\xi\| < \epsilon\}$ is $\Delta^{\ell*}$.
    \end{enumerate}
\end{theorem}

\begin{proof}
    We first prove (i).
    For the first direction, let us assume that $\xi \in E_r$ and let $p \in \beta S$ be an arbitrary idempotent.
    Since $\pi(S)\xi$ is relatively compact, $p-\lim_s\pi(s)\xi$ is well defined in the norm topology, and its value must agree with the value obtained from convergence in the weak topology.
    The desired result now follows from Lemma \ref{UltrafilterLemma}(i).
    For the next direction, let us assume that $\xi\notin E_r$, so that
$\xi=\xi_r+\xi_{\mathrm{aws}}$ for some $\xi_r\in E_r$ and
$0\neq\xi_{\mathrm{aws}}\in E_{\mathrm{aws}}$.
Let $p\in K(\beta S)$ be an idempotent and put
$z:=\xi-\xi_r\neq0$.
Choose $f'\in E'$ such that $\|f'\|=1$ and
$f'(z)=\|z\|$, and set
\[
\epsilon:=\frac12\|z\|>0.
\]
Theorem \ref{ProjectionsViaUltrafilters}(i) tells us that
$p-\lim_s\pi_s\xi=\xi_r$. Hence
\[
p-\lim_s f'(\xi-\pi_s\xi)
=f'(\xi-\xi_r)
=f'(z)
=\|z\|.
\]
It follows that
\[
\{s\in S:\|\xi-\pi_s\xi\|<\epsilon\}\notin p,
\]
which yields the desired result.
    We now prove (ii).
    If $\xi \in E_{\text{aws}}$, then the desired result follows from Lemma \ref{UltrafilterLemma}(iii).
    If $\xi \notin E_{\text{aws}}$, then let $p \in K(\beta S)$ be an idempotent and use Theorem \ref{ProjectionsViaUltrafilters}(i) to see that $p-\lim_s\pi_s\xi \neq 0$.

    Lastly, we prove (iii).
    For the first direction, let us assume that $\xi \in E_r$ and let $p_1,\dots,p_\ell \in S^*$ be such that for $q := p_1^{-1}p_1\dots p_\ell^{-1}p_\ell$ we have $q \in \beta S$.
    Since $\pi(S)\xi$ is relatively compact, $q-\lim_s\pi(s)\xi$ is well defined in the norm topology, and its value must agree with the value obtained from convergence in the weak topology.
    The proof of Lemma \ref{UltrafilterLemma}(ii) can easily be modified to show that $q-\lim_s\pi(s)\xi = \xi$ as desired.
For the next direction, let us assume that $\xi\notin E_r$, and let
$\xi_r$ denote the reversible component of $\xi$.
Put $z:=\xi-\xi_r\neq0$.
Choose $f'\in E'$ such that $\|f'\|=1$ and
$f'(z)=\|z\|$, and set $\epsilon:=\frac12\|z\|>0.$
Let $p_1,\dots,p_\ell\in K(\beta S)$ be arbitrary and put $
q:=p_1^{-1}p_1\cdots p_\ell^{-1}p_\ell.$
By $\ell$ applications of Theorem
\ref{ProjectionsViaUltrafilters}(iii), we have $q-\lim_s\pi_s\xi=\xi_r.$ Hence
\[
q-\lim_s f'(\xi-\pi_s\xi)
=f'(\xi-\xi_r)
=\|z\|.
\]
It follows that
\[
\{s\in S:\|\xi-\pi_s\xi\|<\epsilon\}\notin q,
\]
which yields the desired result.
\end{proof}

\begin{remark}
    When we have a unitary representation of $\mathbb{Z}$ on a Hilbert space, \cite[Theorem 6.2]{bergelson2021iterated} gives a polynomial analogue of Theorem \ref{thm:CharacterizationsViaVisits}(iii) using $\Delta_\ell$ sets (which are special cases of $\Delta^\ell$ sets).
    In light of Theorem \ref{(Weak)EquivalenceForReversiblePart}, it is probable that Theorem \ref{thm:CharacterizationsViaVisits}(iii) has a polynomial generalization for a suitably defined class of polynomial functions over $S$.
    Similarly, \cite[Corollary 6.8]{bergelson2021iterated} gives an analogue of Theorem \ref{thm:CharacterizationsViaVisits}(iii) for weakly mixing measure-preserving systems.
    However, we do not believe that in the more general setting of JdLG-admissible representations on Banach spaces, one can find a polynomial analogue of Theorem \ref{thm:CharacterizationsViaVisits} for the almost weakly stable part.
    Our next example illustrates a crucial difference between weakly mixing vectors in a Hilbert space and almost weakly stable vectors in a Banach space.
\end{remark}

\begin{example}
    Let $\pi$ be the left translation representation of $\mathbb{Z}$ on $E := W(\mathbb{Z})$.
    In this case, the reversible part consists of the space of strongly almost periodic functions, and the almost weakly stable part consists of the functions $f$ for which $|f|$ has $0$ mean.
    Let $S_k \subseteq \mathbb{N}$ denote the set of perfect $k^\text{th}$ powers for some $k \ge 2$.
    Ruppert \cite[Example 11]{OnWeaklyAlmostPeriodicSets} observed that $\mathbbm{1}_{S_k} \in W(\mathbb{Z})$, and it is clear that the mean of $\mathbbm{1}_{S_k}$ is $0$, so $\mathbbm{1}_{S_k}$ is in the almost weakly stable part.
    Let $f' \in E'$ be the evaluation at $0$ map, and observe that $f'(\pi_{n^k}\mathbbm{1}_{S_k}) = 1$ for every $n \in \mathbb{N}$.
    This should be contrasted with the fact that for a unitary operator $U$ on a Hilbert space $\mathcal{H} = \mathcal{H}_c\oplus\mathcal{H}_w$, and vectors $\xi \in \mathcal{H}_w$ and $\eta \in \mathcal{H}$, we have

    \begin{equation}
        \lim_{N\rightarrow\infty}\frac{1}{N}\sum_{n = 1}^N\left|\langle U^{n^k}\xi,\eta\rangle\right| = 0.
    \end{equation}
    See \cite{AlmostWeakPolynomialStabilityOfOperators} for related results and examples.
\end{example}

\normalem
\bibliographystyle{abbrv}
\begin{center}
    \bibliography{references}

@article {kuz-means,
    AUTHOR = {Kuznetsova, Yulia},
     TITLE = {Invariant means on subspaces of quantum weakly almost periodic
              functionals},
   JOURNAL = {J. Anal. Math.},
  FJOURNAL = {Journal d'Analyse Math\'ematique},
    VOLUME = {155},
      YEAR = {2025},
    NUMBER = {1},
     PAGES = {287--316},
      ISSN = {0021-7670,1565-8538},
   MRCLASS = {46L52 (22D25 43A07 43A60 46L07)},
  MRNUMBER = {4905312},
       DOI = {10.1007/s11854-024-0355-y},
       URL = {https://doi.org/10.1007/s11854-024-0355-y},
}

@book {PatersonAmenability,
    AUTHOR = {Paterson, Alan L. T.},
     TITLE = {Amenability},
    SERIES = {Mathematical Surveys and Monographs},
    VOLUME = {29},
 PUBLISHER = {American Mathematical Society, Providence, RI},
      YEAR = {1988},
     PAGES = {xx+452},
      ISBN = {0-8218-1529-6},
   MRCLASS = {43-02 (22-02 43A07 46Lxx)},
  MRNUMBER = {961261},
MRREVIEWER = {C.\ Chou},
       DOI = {10.1090/surv/029},
       URL = {https://doi.org/10.1090/surv/029},
}

@book {AlgebraInTheSCC,
    AUTHOR = {Hindman, Neil and Strauss, Dona},
     TITLE = {Algebra in the {S}tone-\v{C}ech compactification: Theory and applications},
    SERIES = {De Gruyter Textbook},
   EDITION = {Second revised and extended},
 PUBLISHER = {Walter de Gruyter \& Co., Berlin},
      YEAR = {2012},
     PAGES = {xviii+591},
      ISBN = {978-3-11-025623-9},
   MRCLASS = {54-02 (03E05 22A15 54D35 54H99)},
  MRNUMBER = {2893605},
}

@incollection {ERTAnUpdate,
    AUTHOR = {Bergelson, Vitaly},
     TITLE = {Ergodic {R}amsey theory---an update},
 BOOKTITLE = {Ergodic theory of {${\bf Z}^d$} actions ({W}arwick,
              1993--1994)},
    SERIES = {London Math. Soc. Lecture Note Ser.},
    VOLUME = {228},
     PAGES = {1--61},
 PUBLISHER = {Cambridge Univ. Press, Cambridge},
      YEAR = {1996},
   MRCLASS = {28D05 (05A18 05D10 11B25)},
  MRNUMBER = {1411215},
MRREVIEWER = {Karl Petersen},
       DOI = {10.1017/CBO9780511662812.002},
       URL = {https://doi-org.ezproxy.bu.edu/10.1017/CBO9780511662812.002},
}

@article {TheErdosSumsetPaper,
    AUTHOR = {Moreira, Joel and Richter, Florian K. and Robertson, Donald},
     TITLE = {A proof of a sumset conjecture of {E}rd{\"o}s},
   JOURNAL = {Ann. of Math. (2)},
  FJOURNAL = {Annals of Mathematics. Second Series},
    VOLUME = {189},
      YEAR = {2019},
    NUMBER = {2},
     PAGES = {605--652},
      ISSN = {0003-486X},
   MRCLASS = {05D10 (11P70 37A99 46C99)},
  MRNUMBER = {3919363},
MRREVIEWER = {Song Guo},
       DOI = {10.4007/annals.2019.189.2.4},
       URL = {https://doi-org.proxy.lib.ohio-state.edu/10.4007/annals.2019.189.2.4},
}

@article{entangled,
author = {Eisner, Tanja and Kunszenti-Kov{\'a}cs, D{\'a}vid},
 title = {On the pointwise entangled ergodic theorem},
 fjournal = {Journal of Mathematical Analysis and Applications},
 journal = {J. Math. Anal. Appl.},
 issn = {0022-247X},
 volume = {449},
 number = {2},
 pages = {1754--1769},
 year = {2017},
 language = {English},
 doi = {10.1016/j.jmaa.2017.01.011},
 url = {real.mtak.hu/71567/1/1509.05554v4.pdf},
 zbMATH = {6682121},
 Zbl = {1373.37013}
}

@article{Btkai2011Decomposition,
  title={Decomposition of operator semigroups on {W}{*}-algebras},
  author={Andr{\'a}s B{\'a}tkai and Ulrich Groh and D{\'a}vid Kunszenti-Kov{\'a}cs and Marco Schreiber},
  journal={Semigroup Forum},
  year={2011},
  volume={84},
  pages={8-24},
  url={https://api.semanticscholar.org/CorpusID:119115463},
}

@article{MickySohail,
title = {Uniform vector-valued pointwise ergodic theorems for operators},
journal = {Discrete and Continuous Dynamical Systems},
volume = {45},
number = {10},
pages = {3565-3589},
year = {2025},
issn = {1078-0947},
doi = {10.3934/dcds.2025032},
url = {https://www.aimsciences.org/article/id/67d2af282b0ae97dd29c4714},
author = {Micky Barthmann and Sohail Farhangi}
}

@article {EdekoKreidlerHaase,
    AUTHOR = {Edeko, N. and Haase, M. and Kreidler, H.},
     TITLE = {A decomposition theorem for unitary group representations on
              {K}aplansky-{H}ilbert modules and the {F}urstenberg-{Z}immer
              structure theorem},
   JOURNAL = {Anal. Math.},
  FJOURNAL = {Analysis Mathematica},
    VOLUME = {50},
      YEAR = {2024},
    NUMBER = {2},
     PAGES = {377--454},
      ISSN = {0133-3852,1588-273X},
   MRCLASS = {46H25 (06E15 37A15 43A65)},
  MRNUMBER = {4776435},
MRREVIEWER = {R.\ Prakash},
       DOI = {10.1007/s10476-024-00020-1},
       URL = {https://doi.org/10.1007/s10476-024-00020-1},
}

@article {DeLeeuwGlicksberg,
    AUTHOR = {de Leeuw, K. AND Glicksberg, I.},
     TITLE = {Applications of almost periodic compactifications},
   JOURNAL = {Acta Math.},
    VOLUME = {105},
      YEAR = {1961},
     PAGES = {63-97},
  
}

@book {EN,
    AUTHOR = {Engel, Klaus-Jochen and Nagel, Rainer},
     TITLE = {One-parameter semigroups for linear evolution equations},
    SERIES = {Graduate Texts in Mathematics},
    VOLUME = {194},
 PUBLISHER = {Springer-Verlag, New York},
      YEAR = {2000},
     PAGES = {xxii+586},
      ISBN = {0-387-98463-1},
   MRCLASS = {47D06 (34G10 35K90 47N20)},
  MRNUMBER = {1721989},
MRREVIEWER = {Charles\ Batty},
}

@book {OTAoET,
    AUTHOR = {Eisner, Tanja and Farkas, B\'{a}lint and Haase, Markus and Nagel,
              Rainer},
     TITLE = {Operator theoretic aspects of ergodic theory},
    SERIES = {Graduate Texts in Mathematics},
    VOLUME = {272},
 PUBLISHER = {Springer, Cham},
      YEAR = {2015},
     PAGES = {xviii+628},
      ISBN = {978-3-319-16897-5; 978-3-319-16898-2},
   MRCLASS = {47-02 (37-02 37A30 37A55 47A35)},
  MRNUMBER = {3410920},
MRREVIEWER = {Idris Assani},
       DOI = {10.1007/978-3-319-16898-2},
       URL = {https://doi-org.proxy.lib.ohio-state.edu/10.1007/978-3-319-16898-2},
}

@book {UKErgodicThms,
    AUTHOR = {Krengel, Ulrich},
     TITLE = {Ergodic theorems},
    SERIES = {De Gruyter Studies in Mathematics},
    VOLUME = {6},
      NOTE = {With a supplement by Antoine Brunel},
 PUBLISHER = {Walter de Gruyter \& Co., Berlin},
      YEAR = {1985},
     PAGES = {viii+357},
      ISBN = {3-11-008478-3},
   MRCLASS = {28-02 (28Dxx 47A35)},
  MRNUMBER = {797411},
MRREVIEWER = {E. Flytzanis},
       DOI = {10.1515/9783110844641},
       URL = {https://doi-org.proxy.lib.ohio-state.edu/10.1515/9783110844641},
}

@article {RigidityAndNonRecurrence,
    AUTHOR = {Bergelson, V. and del Junco, A. and Lema\'{n}czyk, M. and
              Rosenblatt, J.},
     TITLE = {Rigidity and non-recurrence along sequences},
   JOURNAL = {Ergodic Theory Dynam. Systems},
  FJOURNAL = {Ergodic Theory and Dynamical Systems},
    VOLUME = {34},
      YEAR = {2014},
    NUMBER = {5},
     PAGES = {1464--1502},
      ISSN = {0143-3857},
   MRCLASS = {37A05 (37A25)},
  MRNUMBER = {3255429},
MRREVIEWER = {Lanyu Wang},
       DOI = {10.1017/etds.2013.5},
       URL = {https://doi-org.proxy.lib.ohio-state.edu/10.1017/etds.2013.5},
}

@article {OnModulatedErgodicTheorems,
     author = {Eisner, Tanja and Lin, Michael},
 title = {On modulated ergodic theorems},
 fjournal = {Journal of Nonlinear and Variational Analysis},
 journal = {J. Nonlinear Var. Anal.},
 issn = {2560-6921},
 volume = {2},
 number = {2},
 pages = {131--154},
 year = {2018},
 language = {English},
 doi = {10.23952/jnva.2.2018.2.03},
 zbMATH = {7015091}
   
}

@article{koopman1932dynamical,
  title={Dynamical systems of continuous spectra},
  author={Koopman, Bernard O and Neumann, J v},
  journal={Proceedings of the National Academy of Sciences},
  volume={18},
  number={3},
  pages={255--263},
  year={1932}
}

@article {GodementDecomposition,
    AUTHOR = {Godement, Roger},
     TITLE = {Les fonctions de type positif et la th\'eorie des groupes},
   JOURNAL = {Trans. Amer. Math. Soc.},
  FJOURNAL = {Transactions of the American Mathematical Society},
    VOLUME = {63},
      YEAR = {1948},
     PAGES = {1--84},
      ISSN = {0002-9947,1088-6850},
   MRCLASS = {20.0X},
  MRNUMBER = {23243},
MRREVIEWER = {I.\ E.\ Segal},
       DOI = {10.2307/1990636},
       URL = {https://doi.org/10.2307/1990636},
}

@book {TempelmanBook,
    AUTHOR = {Tempelman, Arkady},
     TITLE = {Ergodic theorems for group actions},
    SERIES = {Mathematics and its Applications},
    VOLUME = {78},
      NOTE = {Informational and thermodynamical aspects,
              Translated and revised from the 1986 Russian original},
 PUBLISHER = {Kluwer Academic Publishers Group, Dordrecht},
      YEAR = {1992},
     PAGES = {xviii+399},
      ISBN = {0-7923-1717-3},
   MRCLASS = {22D40 (28D15 47A35 47D03 60G60)},
  MRNUMBER = {1172319},
MRREVIEWER = {Arlan\ Ramsay},
       DOI = {10.1007/978-94-017-1460-0},
       URL = {https://doi.org/10.1007/978-94-017-1460-0},
}

@article{dye1965ergodic,
  title={On the ergodic mixing theorem},
  author={Dye, HA},
  journal={Transactions of the American Mathematical Society},
  volume={118},
  pages={123--130},
  year={1965},
  publisher={JSTOR}
}

@article {JacobsJdLG,
    AUTHOR = {Jacobs, Konrad},
     TITLE = {Ergodentheorie und fastperiodische {F}unktionen auf
              {H}albgruppen},
   JOURNAL = {Math. Z.},
  FJOURNAL = {Mathematische Zeitschrift},
    VOLUME = {64},
      YEAR = {1956},
     PAGES = {298--338},
      ISSN = {0025-5874,1432-1823},
   MRCLASS = {46.2X},
  MRNUMBER = {77092},
MRREVIEWER = {M.\ M.\ Day},
       DOI = {10.1007/BF01166575},
       URL = {https://doi.org/10.1007/BF01166575},
}

@misc{RNFixedPoint,
 author = {Ryll-Nardzewski, Czes{\l}aw},
 title = {On fixed points of semigroups of endomorphisms of linear spaces},
 year = {1967},
 language = {English},
 howpublished = {Proc. 5th {Berkeley} {Symp}. {Math}. {Stat}. {Probab}., {Univ}. {Calif}. 1965/66, 2, {Part} 1, 55-61},
 zbMATH = {3301909},
 Zbl = {0189.17501}
}

@book{assani2003wiener,
  title={Wiener-Wintner ergodic theorems},
  author={Assani, Idris},
  year={2003},
  publisher={World Scientific Publishing Company}
}

@book{AnalysisOnSemigroups,
 author = {Berglund, John F. and Junghenn, Hugo D. and Milnes, Paul},
 title = {Analysis on semigroups: function spaces, compactifications, representations},
 fseries = {Canadian Mathematical Society Series of Monographs and Advanced Texts},
 series = {Can. Math. Soc. Ser. Monogr. Adv. Texts},
 isbn = {0-471-61208-1},
 year = {1989},
 publisher = {New York etc.: John Wiley \& Sons},
 language = {English},
 zbMATH = {44681},
 Zbl = {0727.22001}
}

@article{day1950means,
  title={Means for the bounded functions and ergodicity of the bounded representations of semi-groups},
  author={Day, Mahlon M},
  journal={Transactions of the American Mathematical Society},
  volume={69},
  number={2},
  pages={276--291},
  year={1950},
  publisher={JSTOR}
}

@article{baker1976stone,
  author = {Baker, J. W. and Butcher, R. J.},
 title = {The {Stone}-{Cech} compactification of a topological semigroup},
 fjournal = {Mathematical Proceedings of the Cambridge Philosophical Society},
 journal = {Math. Proc. Camb. Philos. Soc.},
 issn = {0305-0041},
 volume = {80},
 pages = {103--107},
 year = {1976},
 language = {English},
 doi = {10.1017/S0305004100052725},
 zbMATH = {3515741},
 Zbl = {0329.22003}
}

@article{BMIPUltrafiltersAndSzemerediForGPolynomials,
 author = {Bergelson, Vitaly and McCutcheon, Randall},
 title = {Idempotent ultrafilters, multipleweak mixing and {Szemer{\'e}di}'s theorem for generalized polynomials},
 fjournal = {Journal d'Analyse Math{\'e}matique},
 journal = {J. Anal. Math.},
 issn = {0021-7670},
 volume = {111},
 pages = {77--130},
 year = {2010},
 language = {English},
 doi = {10.1007/s11854-010-0013-4},
 zbMATH = {5864337},
 Zbl = {1220.37003}
}

@article{ehrenpreis1955uniformly,
  title={Uniformly bounded representations of groups},
  author={Ehrenpreis, Leon and Mautner, Friederich I},
  journal={Proceedings of the National Academy of Sciences},
  volume={41},
  number={4},
  pages={231--233},
  year={1955}
}

@inproceedings{glasscock2025folner,
  title={F{\o}lner, {B}anach, and translation density are equal and other new results about density in left amenable semigroups: D. {G}lasscock et al.},
  author={Glasscock, Daniel and Hindman, Neil and Strauss, Dona},
  booktitle={Semigroup Forum},
  pages={1--23},
  year={2025},
  organization={Springer}
}

@article{beiglbock2009solvability,
  title={Solvability of Rado systems in D-sets},
  author={Beiglb{\"o}ck, Mathias and Bergelson, Vitaly and Downarowicz, Tomasz and Fish, Alexander},
  journal={Topology and its Applications},
  volume={156},
  number={16},
  pages={2565--2571},
  year={2009},
  publisher={Elsevier}
}

@article{bergelson2016polynomial,
  title={Polynomial recurrence with large intersection over countable fields},
  author={Bergelson, Vitaly and Robertson, Donald},
  journal={Israel Journal of Mathematics},
  volume={214},
  number={1},
  pages={109--120},
  year={2016},
  publisher={Springer}
}

@article{campbell2018d,
  author = {Campbell, James T. and McCutcheon, Randall},
 title = {D sets and {IP} rich sets in countable, cancellative abelian semigroups},
 fjournal = {Semigroup Forum},
 journal = {Semigroup Forum},
 issn = {0037-1912},
 volume = {96},
 number = {1},
 pages = {49--68},
 year = {2018},
 language = {English},
 doi = {10.1007/s00233-017-9854-9},
 zbMATH = {6888055},
 Zbl = {1456.20063}
}

@article{mccutcheon2014d,
  title={${D}$ sets and a {S}{\'a}rk{\"o}zy theorem for countable fields},
  author={McCutcheon, Randall and Windsor, Alistair},
  journal={Israel Journal of Mathematics},
  volume={201},
  number={1},
  pages={123--146},
  year={2014},
  publisher={Springer}
}

@article{banakh2022automatic,
  title={Automatic continuity of measurable homomorphisms on Cech-complete topological groups},
  author={Banakh, Taras},
  journal={arXiv preprint arXiv:2206.02481},
  year={2022}
}

@article{jones1976ergodic,
  title={Ergodic theorems of weak mixing type},
  author={Jones, Lee K and Lin, Michael},
  journal={Proceedings of the American Mathematical Society},
  volume={57},
  number={1},
  pages={50--52},
  year={1976}
}

@article{eisner2007weakly,
  title={Weakly and almost weakly stable ${C}_0$-semigroups},
  author={Eisner, Tanja and Farkas, Balint and Nagel, Rainer and Sereny, Andras},
  journal={International Journal of Dynamical Systems and Differential Equations},
  volume={1},
  number={1},
  pages={44--57},
  year={2007},
  publisher={Inderscience Publishers}
}

@article{bergelson2021iterated,
  title={Iterated differences sets, Diophantine approximations and applications},
  author={Bergelson, Vitaly and Zelada, Rigoberto},
  journal={Journal of Combinatorial Theory, Series A},
  volume={184},
  pages={105520},
  year={2021},
  publisher={Elsevier}
}

@article{Rudin1959,
  author    = {Walter Rudin},
  title     = {Weak Almost Periodic Functions and Fourier-Stieltjes Transforms},
  journal   = {Duke Mathematical Journal},
  volume    = {26},
  year      = {1959},
  pages     = {215--220},
  doi       = {10.1215/S0012-7094-59-02621-4},
  mrnumber  = {0107265}
}

@article{KornfeldLin2000,
  author  = {Isaac Kornfeld and Michael Lin},
  title   = {Weak Almost Periodicity of {$L_1$} Contractions and Coboundaries of Non-Singular Transformations},
  journal = {Studia Mathematica},
  volume  = {138},
  number  = {3},
  pages   = {225--240},
  year    = {2000},
  doi     = {10.4064/sm-138-3-225-240}
}

@book{DunklRamirez1971,
  author    = {Charles F. Dunkl and Donald E. Ramirez},
  title     = {Topics in Harmonic Analysis},
  publisher = {Appleton-Century-Crofts},
  address   = {New York},
  year      = {1971}
}

@article{klawe1980dimensions,
  title={Dimensions of the sets of invariant means of semigroups},
  author={Klawe, Maria M},
  journal={Illinois Journal of Mathematics},
  volume={24},
  number={2},
  pages={233--243},
  year={1980},
  publisher={Duke University Press}
}

@article{bergelson2007central,
  title={Central sets and a non-commutative Roth theorem},
  author={Bergelson, Vitaly and McCutcheon, Randall},
  journal={American journal of mathematics},
  volume={129},
  number={5},
  pages={1251--1275},
  year={2007},
  publisher={Johns Hopkins University Press}
}

@misc{burckel1970weakly,
 author = {Burckel, R. B.},
 title = {Weakly almost periodic functions on semigroups},
 year = {1970},
 language = {English},
 howpublished = {Notes on {Mathematics} and its {Applications}. {New} {York}-{London}-{Paris}: {Gordon} and {Breach} {Science} {Publishers}. ix, 118 p.},
 zbMATH = {3306233},
 Zbl = {0192.48602}
}

@incollection{MinimalIdempotentsAndERT,
 author = {Bergelson, Vitaly},
 title = {Minimal idempotents and ergodic {Ramsey} theory},
 booktitle = {Topics in dynamics and ergodic theory. Survey papers and mini-courses presented at the international conference and US-Ukrainian workshop on dynamical systems and ergodic theory, Katsiveli, Ukraine, August 21--30, 2000},
 isbn = {0-521-53365-1},
 pages = {8--39},
 year = {2003},
 publisher = {Cambridge: Cambridge University Press},
 language = {English},
 zbMATH = {2070330},
 Zbl = {1039.05063}
}

@article{QuotientSetsAndDensityRecurrentSets,
 author = {Bergelson, Vitaly and Hindman, Neil},
 title = {Quotient sets and density recurrent sets},
 fjournal = {Transactions of the American Mathematical Society},
 journal = {Trans. Am. Math. Soc.},
 issn = {0002-9947},
 volume = {364},
 number = {9},
 pages = {4495--4531},
 year = {2012},
 language = {English},
 doi = {10.1090/S0002-9947-2012-05417-1},
 zbMATH = {6191419},
 Zbl = {1273.22002}
}

@book {AlgebraicTheoryOfSemigroupsVol1,
    AUTHOR = {Clifford, A. H. and Preston, G. B.},
     TITLE = {The algebraic theory of semigroups. {V}ol. {I}},
    SERIES = {Mathematical Surveys, No. 7},
 PUBLISHER = {American Mathematical Society, Providence, RI},
      YEAR = {1961},
     PAGES = {xv+224},
   MRCLASS = {20.92},
  MRNUMBER = {132791},
MRREVIEWER = {\v{S}t. Schwarz},
}

@book {AmenabilityByPaterson,
    AUTHOR = {Paterson, Alan L. T.},
     TITLE = {Amenability},
    SERIES = {Mathematical Surveys and Monographs},
    VOLUME = {29},
 PUBLISHER = {American Mathematical Society, Providence, RI},
      YEAR = {1988},
     PAGES = {xx+452},
      ISBN = {0-8218-1529-6},
   MRCLASS = {43-02 (22-02 43A07 46Lxx)},
  MRNUMBER = {961261},
MRREVIEWER = {C. Chou},
       DOI = {10.1090/surv/029},
       URL = {https://doi-1org-102eb17k00086.han.amu.edu.pl/10.1090/surv/029},
}

@book{RandomPaperForOreLemma,
 author = {Davidson, Kenneth R. and Fuller, Adam H. and Kakariadis, Evgenios T. A.},
 title = {Semicrossed products of operator algebras by semigroups},
 fseries = {Memoirs of the American Mathematical Society},
 series = {Mem. Am. Math. Soc.},
 issn = {0065-9266},
 volume = {1168},
 isbn = {978-1-4704-2309-4; 978-1-4704-3697-1},
 year = {2017},
 publisher = {Providence, RI: American Mathematical Society (AMS)},
 language = {English},
 doi = {10.1090/memo/1168},
 zbMATH = {6724255},
 Zbl = {1385.47027}
}

@article{OnWeaklyAlmostPeriodicSets,
 author = {Ruppert, W. A. F.},
 title = {On weakly almost periodic sets},
 fjournal = {Semigroup Forum},
 journal = {Semigroup Forum},
 issn = {0037-1912},
 volume = {32},
 pages = {267--281},
 year = {1985},
 language = {English},
 doi = {10.1007/BF02575544},
 url = {https://eudml.org/doc/134786},
 zbMATH = {3908716},
 Zbl = {0569.22001}
}

@article{AlmostWeakPolynomialStabilityOfOperators,
 author = {Kunszenti-Kov{\'a}cs, D{\'a}vid},
 title = {Almost weak polynomial stability of operators},
 fjournal = {Houston Journal of Mathematics},
 journal = {Houston J. Math.},
 issn = {0362-1588},
 volume = {41},
 number = {3},
 pages = {901--913},
 year = {2015},
 language = {English},
 url = {math.uh.edu/~hjm/restricted/pdf41(3)/10kunszenti.pdf},
 zbMATH = {6514335},
 Zbl = {1329.47016}
}

@book{TanjaEisnerStabilityBook,
 author = {Eisner, Tanja},
 title = {Stability of operators and operator semigroups},
 fseries = {Operator Theory: Advances and Applications},
 series = {Oper. Theory: Adv. Appl.},
 issn = {0255-0156},
 volume = {209},
 isbn = {978-3-0346-0194-8},
 year = {2010},
 publisher = {Basel: Birkh{\"a}user},
 language = {English},
 zbMATH = {5625330},
 Zbl = {1205.47002}
}

@article{CCCEberleinCompactaAreMetrizable,
 author = {Argyros, S. and Negrepontis, S.},
 title = {On weakly {K}-countably determined spaces of continuous functions},
 fjournal = {Proceedings of the American Mathematical Society},
 journal = {Proc. Am. Math. Soc.},
 issn = {0002-9939},
 volume = {87},
 pages = {731--736},
 year = {1983},
 language = {English},
 doi = {10.2307/2043370},
 zbMATH = {3804416},
 Zbl = {0509.54020}
}
\end{center}

\end{document}